\documentclass[12pt,a4paper]{article}

\usepackage[dvipdfmx]{color}

\usepackage[pdftex]{graphicx}

\usepackage{authblk}

\usepackage{amsmath}
\usepackage{amsthm}
\usepackage{ascmac}
\usepackage{bm}
\usepackage{cancel}
\usepackage{comment}
\usepackage{fancyhdr}
\usepackage{cases}
\usepackage{multicol}
\usepackage{braket}
\usepackage{geometry}
\theoremstyle{definition}
\newtheorem{theorem}{Theorem}
\newtheorem*{theorem*}{Theorem}
\newtheorem{definition}{Definition}
\newtheorem*{definition*}{Definition}
\newtheorem{remark}{Remark}
\newtheorem*{remark*}{Remark}
\newtheorem{proposition}{Proposition}
\newtheorem*{proposition*}{Proposition}
\newtheorem*{thm}{Theorem}

\title{Differential and Integral Calculus of Sequence}
\author{Yusuke Imai \thanks{CONTACT: 93imaiyusuke@gmail.com}}
\affil{Graduate School of Engineering Science, Osaka University, Toyonaka, Osaka 560-8531, Japan}
\date{\today}

\begin{document}
\maketitle

\begin{abstract}
We create a sequence version of calculus. First, we define equivalence, some fundamental operations, differential, and integral for sequences. Then, we propose sequence versions of identity function, power function, exponential function, hyperbolic function, trigonometric function, and also find sequence versions of the Maclaurin series for them. The sequence versions of exponential function involve divergent series including Grandi's series. By using this framework, we find a sequence version of the binomial theorem and Euler's identity. In addition, we design new formalisms of Fibonacci sequence and its generalizations. Last, we propose a sequence dual of factorial and Bell number, and find sequence dual of modular property of factorial concerning prime number (Wilson's theorem) and of Bell number concerning prime number. 
\end{abstract}
 
 \tableofcontents

\newpage
\section{Introduction} 
A sequence is an ordered list of numbers and is a fundamental mathematical concept that has been investigated since ancient times. The following sequence would be the most basic sequence.
  \begin{align}
     1, 1, 1, 1, 1, 1, \cdots.
  \end{align}
This sequence has the following three properties. First, all the numbers in the above sequence are 1. Second, all the differences between the two consecutive terms are 0. By this property, we can see that the above sequence is an arithmetic sequence. Third, all the ratios between the two consecutive terms are 1. By this property, we can see that the above sequence is a geometric sequence. The following sequence is also a fundamental sequence that is an arithmetic sequence whose common difference is 1 and the initial term is given by 0. 
  \begin{align}
    0, 1, 2, 3, 4, 5, \cdots.
  \end{align}
The following sequence is also a fundamental sequence that is a geometric sequence whose common ratio is 2 and the initial term is given by 1. 
  \begin{align}
    1, 2, 4, 8, 16, 32, \cdots.
  \end{align}
Each term of arithmetic sequences or geometric sequences is uniquely determined by the previous term except the first term.

There exist sequences, each of whose term is determined by the previous two terms, except the first two terms. The most famous one would be the sequence each of whose term is given by sum of the previous two terms, except the first two terms and that starts with 0 and 1, i.e. the Fibonacci sequence \cite{Fibonacci:history}.
  \begin{align}
    0, 1, 1, 2, 3, 5, \cdots. 
  \end{align}
The most characteristic feature of the Fibonacci sequence is that the ratio of consecutive two terms in the Fibonacci sequence approach to the golden ratio. The Fibonacci sequence was proposed in 1202 by Leonardo of Pisa \cite{Fibonacci:history}, and it has been discovered that it characterizes many things in nature, such as bee populations, spiral patterns of poplar, willow, pear trees, beech, hazel, cherry, apple, elm, lime, and almond \cite{Fibonacci:history2}.

\begin{figure}[h]
 \centering
 \includegraphics[width=14.0cm]{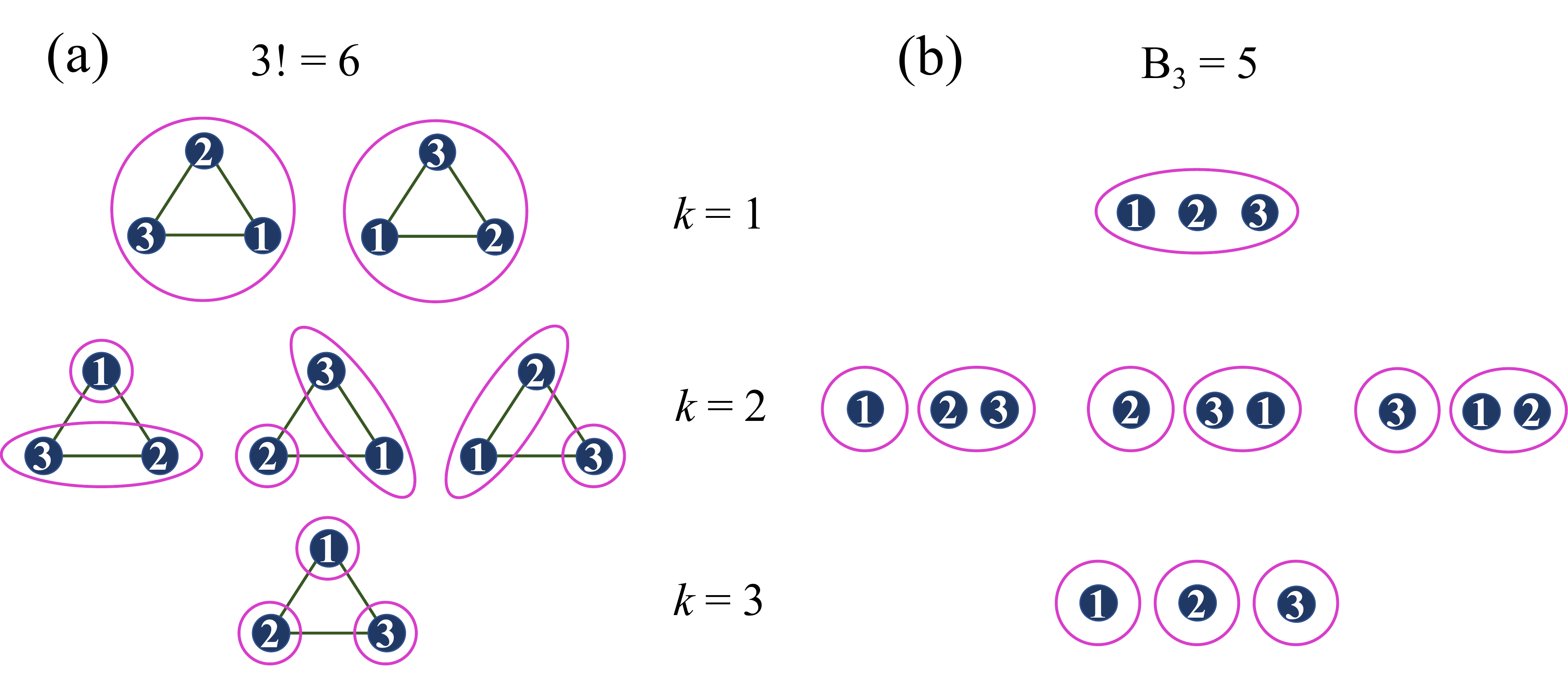}
 \caption{Diagrammatic representation of (a) factorial and the unsigned Stirling number of the first kind and (b) Bell number and Stirling number of the second kind for $n=3$. The left diagram shows permutations of vertices of a triangle for the arrangement of the bottom figure. The permutations are represented by the circles, and the number of the circles corresponds to $k$. The right diagram shows grouping of three numbers. The groups are represented by the circles, and the number of the circles corresponds to $k$.}
 \label{factorial_Bell}
\end{figure}

In addition, there exist sequences, each of whose term is determined by its index and the previous terms, except the first term. For example, each term of the following sequence is given by $n$ times the previous term, except the first term if the first index is set to 0.
  \begin{align}
    1, 1, 2, 6, 24, 120, \cdots.
  \end{align}
The $n$th term in the above sequence equals to factorial of $n$ that represents the number of permutations for $n$ elements. The ($n$, $k$)th unsigned Stirling numbers of the first kind represents the number of permutations for $n$ elements with $k$ disjoint cyclic permutations. Then, sum of the ($n$, $k$)th unsigned Stirling number of the first kind from $k=0$ to $n$ equals to the $n$th factorial \cite{Stirling} (Fig.~\ref{factorial_Bell}(a)).

The following sequence is called Bell numbers and is also a sequence that has combinatorial meaning and can be considered as a dual of factorial in terms of the Stirling numbers \cite{Bell, Bell2, Bell3}. 
  \begin{align}
    1, 1, 2, 5, 15, 52, \cdots.
  \end{align}
The $n$th Bell number represents the number of ways to group $n$ elements. The ($n$, $k$)th Stirling number of the second kind represents the number of ways to group $n$ elements by $k$ disjoint groups. Then, sum of the ($n$, $k$)th Stirling number of the second kind from $k=0$ to $n$ equals to the $n$th Bell number. While the ($n$, $k$)th unsigned Stirling numbers of the first kind can also be defined as coefficients in a linear combination of $x^k$ for $k=0, \cdots, n$ that equals to the rising factorial, $x(x+1)\cdots(x+n-1)$, the ($n$, $k$)th Stirling number of the second kind can also be defined as coefficients in a linear combination of the falling factorial $x(x-1)\cdots(x-k+1)$ for $k=0, \cdots, n$ that equals to $x^n$ \cite{Stirling} (Fig.~\ref{factorial_Bell}(b)). Furthermore, both the factorials and the Bell numbers have the modular property concerning prime numbers. The following theorem is called Wilson's theorem \cite{Wilson, Wilson2, Wilson3} and states that the $(p-1)$th factorial equals to $-1$ modulo $p$ if and only if $p$ is a prime number. 
\begin{thm}
  \begin{align}
    (p-1) ! \equiv -1 \ ({\rm mod} \ p). 
  \label{intro:Wilson}
  \end{align}
\end{thm}
Also, sum of the $n$th Bell number and $(n+1)$th Bell number equals to $(n+p)$th Bell number modulo $p$ if $p$ is a prime number \cite{Bell:prime}.
\begin{thm}
  \begin{align}
    B_{n+p} \equiv B_n + B_{n+1} \ ({\rm mod} \ p).
  \label{intro:Bell}
  \end{align}
\end{thm}

\begin{figure}[h]
 \centering
 \includegraphics[width=8.0cm]{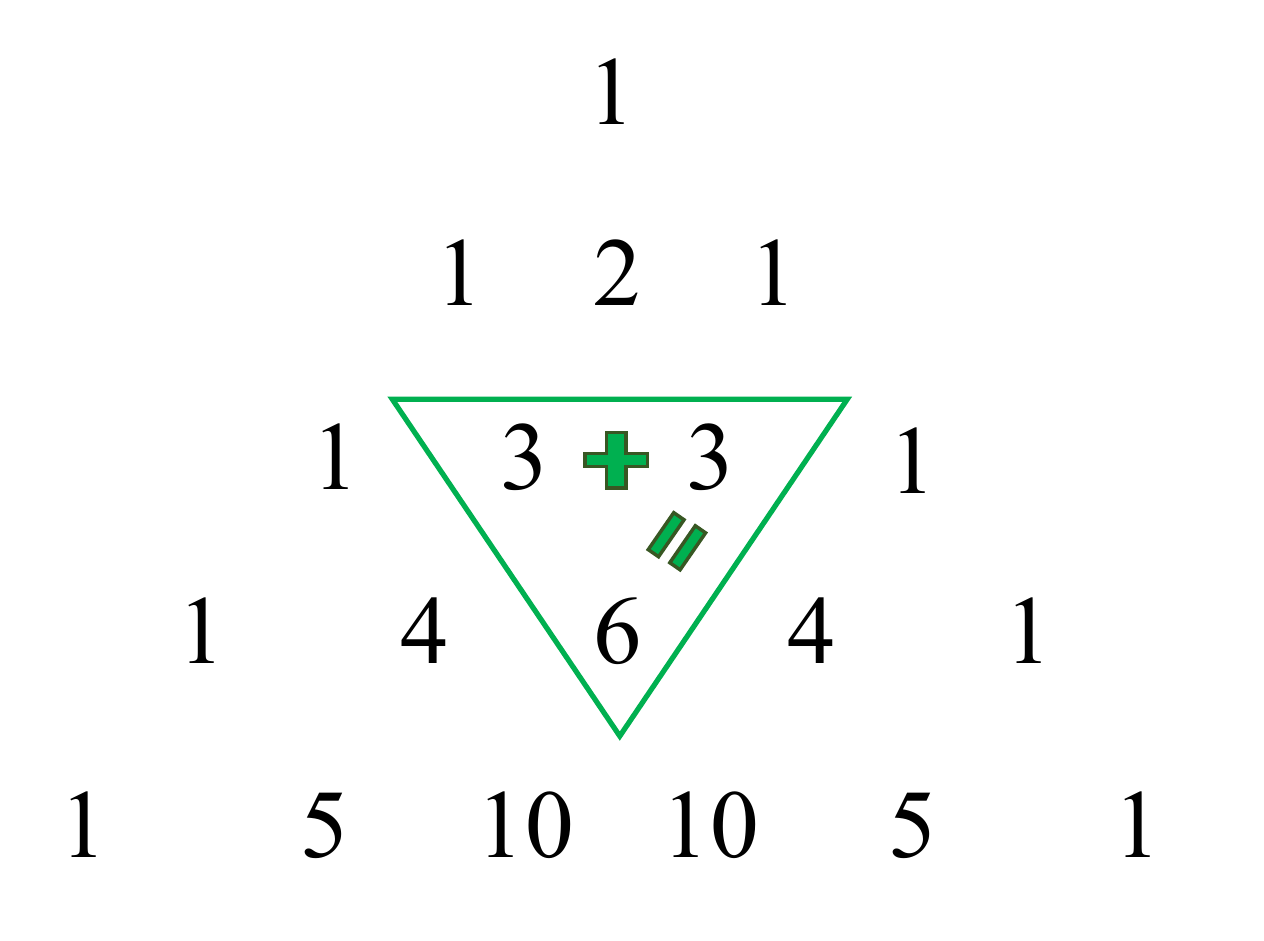}
 \caption{Pascal's triangle.}
 \label{Pascal}
\end{figure}

On the other hand, calculus is a fundamental mathematical branch that deals with change or accumulation of a quantity \cite{calculus}. The following real function would be the most basic function.
  \begin{align}
    f(x) = 1.
  \end{align}
This function has the following three properties. First, all the values the above function takes are 1. Second, $f(x+1) - f(x) = 0$ for all $x$ in real numbers. Third, $f(x+1)/f(x) = 1$ for all $x$ in real numbers. The following real function is also a fundamental real function that is called identity function and that satisfies $f(x+1) - f(x) = 1 $.
  \begin{align}
    f(x) = x.
  \end{align}
The following theorem is called binomial theorem \cite{Pascal:binomial}.
\begin{thm}
  \begin{align}
    (1 + x)^n = \sum_{k = 0}^n {}_nC_k \, x^k. 
  \label{binomial:intro}
  \end{align}
\end{thm}
If $x=1$, the binomial theorem can be visualized by Pascal's triangle (Fig.~\ref{Pascal}).

The following real function is a fundamental real function that is called exponential function and that satisfies $f(x+1)/f(x) = e$ where $e$ is the Napier's constant.
  \begin{align}
    f(x) = e^x.
  \end{align}
The exponential function is often used to represent an explosive increase of quantity such as the number of infected persons \cite{exponential:infection}. The exponential function is also characterized by invariance for differential as follows.
  \begin{align}
    \frac{d}{dx} e^x = e^x.
  \label{ex_differential}
  \end{align}

The hyperbolic functions, $\sinh x$ and $\cosh x$, are determined by the exponential function as follows.
  \begin{align}
    \cosh x &= \frac{e^x + e^{-x}}{2},
  \\
    \sinh x &= \frac{e^x - e^{-x}}{2}.
  \end{align}
Perhaps the most famous application of hyperbolic functions is the four-dimensional rotation in special relativity \cite{hyperbolic:relativity}. The hyperbolic functions satisfy the following property for differential.
  \begin{align}
    \frac{d }{dx} \cosh x &= \sinh x,
  \\
    \frac{d }{dx} \sinh x &= \cosh x.
  \end{align}
Also, the trigonometric functions, $\cos x$ and $\sin x$, are determined by the exponential function as follows.
  \begin{align}
    \cos x &= \frac{e^{ix} + e^{-ix}}{2},
  \\
    \sin x &= \frac{e^{ix} - e^{-ix}}{2i}.
  \end{align}
They are used to represent rotation in Euclidean space \cite{trigonometric:rotation}, and $\cos x$ and $\sin x$ have the following periodic property.
  \begin{align}
    \cos (x+2\pi) &= \cos x,
  \\
    \sin (x+2\pi) &= \sin x.
  \end{align}
Also, the trigonometric functions satisfy the following property for differential.
  \begin{align}
    \frac{d }{dx} \cos x &= -\sin x,
  \\
    \frac{d }{dx} \sin x &= \cos x.
  \end{align}
One can see that the differential rule for the trigonometric functions is similar to one of the hyperbolic functions. The following relations also indicate the similarity between the trigonometric functions and the hyperbolic functions.
  \begin{align}
    &\cos^2 x + \sin^2x = 1,
  \label{intro:tri:1}
  \\
    &\cosh^2 x - \sinh^2x = 1,
  \label{intro:hyp:1}
  \\
    &\cos^2 x - \sin^2x = \cos 2x,
  \label{intro:tri:2}
  \\
    &\cosh^2 x + \sinh^2x = \cosh 2x,
  \label{intro:hyp:2}
  \\
    & 2 \cos x \sin x = \sin 2x,
  \label{intro:tri:3}
  \\
    & 2\cosh x \sinh x = \sinh 2x.
  \label{intro:hyp:3}
  \end{align}
The exponential function, the hyperbolic functions, and the trigonometric functions have the Maclaurin series, i.e. can be represented by the series of $1, x, x^2, \cdots$ with coefficients.
  \begin{align}
    e^x &= \sum_{n=0}^{\infty} \frac{x^n}{n!},
  \label{Maclaurin:ex}
  \\
    \cosh x &= \sum_{n=0}^{\infty} \frac{x^{2n}}{(2n)!},
  \label{Maclaurin:coshx}
  \\
    \sinh x &= \sum_{n=0}^{\infty} \frac{x^{2n+1}}{(2n+1)!},
  \label{Maclaurin:sinhx}
  \\
    \cos x &= \sum_{n=0}^{\infty} \frac{(-1)^n x^{2n}}{(2n)!},
  \label{Maclaurin:cosx}
  \\
    \sin x &= \sum_{n=0}^{\infty} \frac{(-1)^n x^{2n+1}}{(2n+1)!}.
  \label{Maclaurin:sinx}
  \end{align}
By differentiating the above Maclaurin series by term, one can show that they satisfy the differential rules shown above.

\begin{figure}[h]
 \centering
 \includegraphics[width=15cm]{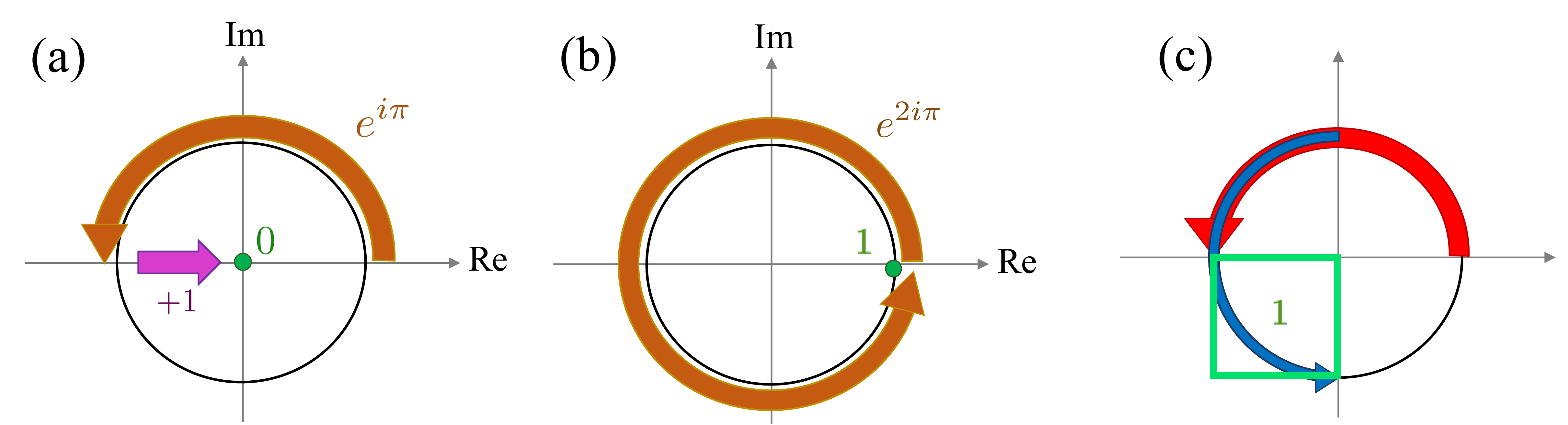}
 \caption{Diagrammatic representation of (a) Euler's identity ($e^{i\pi}+1=0$) and (b) $e^{2i\pi}=1$, and (c) geometric representation of $e^{i\pi}\times e^{i\pi}=1$ (area of the square equals to 1).}
 \label{Euler}
\end{figure}

Also, the exponential function relates to the trigonometric functions as follows.
  \begin{align}
    e^{ix} = \cos x + i \sin x.
  \label{intro:ex:tri}
  \end{align}
The following equality is called Euler's identity \cite{Euler:identity} and would be the most fundamental and beautiful theorem in mathematics because it connects various mathematical concepts, the $e$ (Napier's constant), $0$ (additive identity), $1$ (multiplicative identity), $i$ (imaginary unit), $\pi$ (circle ratio).
\begin{thm}
  \begin{align}
    e^{i \pi} + 1 = 0.
  \label{Euler:intro}
  \end{align}
\end{thm}
Figure \ref{Euler}(a) shows Euler's identity ($e^{i\pi}+1=0$) and Fig.~\ref{Euler}(b) shows $e^{2i\pi}=1$, and Fig.~\ref{Euler}(c) shows geometric representation of $e^{i\pi}\times e^{i\pi}=1$.

In addition, we introduce a function relating to the factorial and a function relating to the Bell number. The following function is called gamma function and can be defined in complex numbers whose real part is positive \cite{gamma}.
  \begin{align}
    \Gamma (z) = \int_0^{\infty} t^{z-1} e^{-t} dt. 
  \label{gamma}
  \end{align}
The gamma function satisfies the following relation
  \begin{align}
    \Gamma(n+1) = n!,
  \end{align}
i.e. the gamma function can be regarded as an extension of the factorial.

Also, the exponential function of the exponential function divided by $e$ is represented as follows \cite{Bell}.
  \begin{align}
     e^{e^x}/e = \sum_{n=0}^{\infty} \frac{B_n x^n}{n!}.
  \label{eex}
  \end{align}
Then, the above function is called the exponential generating function of $B_n$. Note that the exponential generating function of $1,1,1,\cdots$ is given by the exponential function.

This paper investigates relations among the sequences and functions shown above and their generalizations by creating a sequence version of calculus. In Ref.~\cite{Diff:Int:seq}, differential and integral of a sequence were proposed and they are applied to solve discrete differential equations. In this paper, by inheriting the concepts proposed in Ref.~\cite{Diff:Int:seq} with some retouching, we build a sequence version of calculus. Because of the difference between the discreteness and the continuity, it is difficult to construct a sequence version of calculus that is perfectly analogous to usual calculus. However, we will see there exist many beautiful sequence-analogies of calculus.

While writing this paper, we found a note \cite{Yuki:seq} that tried to find a discrete version of some elementary functions: exponential function, logarithmic function, and trigonometric function. The note concludes with the sentences after dealing with the discrete version of the trigonometric functions ($f(x)$ and $g(x)$); ``Looking at this, the author wonders if $f(x)$ and $g(x)$ are really the most appropriate discrete version of the trigonometric functions. Could not a more appropriate function be found? What do the readers think?". Because, in this paper, we propose three kinds of sequence version of the trigonometric functions and one of them corresponds to $f(x)$ and $g(x)$ (Eqs.~(\ref{Yuki1}) and (\ref{Yuki2})) and the three sequence versions of the trigonometric functions are entangled deeply each other to result in various formulae that can be regarded as sequence version of formulae of the trigonometric functions, we wish this paper can be one of the answers to that question.

In Sec.~\ref{Sec:definition}, we define some concepts relating to sequence. First, we show how to represent a sequence and define equivalence between two sequences. Next, we define fundamental operations including shift operator, sum, subtraction, multiplication, division, scalar addition and multiplication, inverse, and insertion. Then, we define two kinds of differential and integral: left differential, right differential, left integral, right integral, and reveal the relations among them.

In Sec.~\ref{Sec:fundamental}, we define the sequence version of fundamental real functions such as additive identity ($f(x)=0$), multiplicative identity ($f(x)=1$), and generally constant function ($f(x)=a$) where $a$ is a real number. We also introduce the general notation that represents the sequence version of a real function.

In Sec.~\ref{Sec:power}, we define the sequence version of the identity function ($f(x)=x$), and generally power function ($f(x)=x^a$) where $a$ is a real number. Also, we show the law of exponent holds for the sequence version of the power function. Then, we obtain a direct analogy of the binomial theorem (Eq.~(\ref{binomial:intro})). This can be proved by using the usual binomial theorem.

In Sec.~\ref{Sec:exponential}, we first define the sequence version of exponential function by focusing on the property of the exponential function for differential (Eq.~(\ref{ex_differential})). Then, we obtain the sequence versions of the Maclaurin series for the exponential function, some of which relates to divergent series including Grandi’s series. Also, we find the right/left integral of the sequence version of $x^n$ relates to $x^{n+1}$ through the Eulerian number. Next, we define two kinds of sequence version of the negative exponential function depending on whether the left differential or right differential is used. Last, we define two kinds of sequence versions of the general exponential function whose exponent is given by a general real number $a$.

In Sec.~\ref{Sec:hyperbolic}, we define two kinds of sequence version of the hyperbolic functions by using the sequence version of the exponential function and the two kinds of the negative exponential function defined in Sec.~\ref{Sec:exponential}. Each of them satisfies the sequence version of Eq.~(\ref{intro:hyp:1}), Eq.~(\ref{intro:hyp:2}), and Eq.~(\ref{intro:hyp:3}).

In Sec.~\ref{Sec:trigonometric}, we first define two kinds of sequence version of the sine function and cosine function by using Eq.~(\ref{intro:ex:tri}) and the two kinds of sequence version of the general exponential function defined in Sec.~\ref{Sec:exponential}. Each of them satisfies the sequence version of Eq.~(\ref{intro:tri:1}), Eq.~(\ref{intro:tri:2}), and Eq.~(\ref{intro:tri:3}). Then, we obtain sequence version of Euler's identity (Eq.~(\ref{Euler:intro})). Also, by combining the two sequence versions of the sine (cosine) function, we obtain another sequence version of the sine (cosine) function whose periodicity is 8.

In Sec.~\ref{Sec:fibonacci}, we propose the sequence version of the Maclaurin series for the Fibonacci sequence, $(P, \, Q)$-Fibonnaci sequence including Pell number ($P=2$, $Q=1$) and Jacobsthal number ($P=1$, $Q=2$), and $k$-bonacci sequence. We find similarities between those generalizations of the Fibonacci sequence and the sequence version of exponential function from the viewpoint of the sequence version of the Maclaurin series.

In Sec.~\ref{Sec:factorial:Bell}, we introduce sequence dual of factorial and of Bell numbers by focusing on the functions shown in Eq.~(\ref{gamma}) and Eq.~(\ref{eex}), and we find that sequence dual of factorial (Bell numbers) is characterized by the Stirling transform of the second (first) kind for the Bell numbers (factorial). In addition, we find that the sequence dual of the factorial (Bell numbers) inherits partially the modular property of the factorial (Bell numbers) for the prime numbers shown in Eq.~(\ref{intro:Wilson}) (Eq.~(\ref{intro:Bell})).

Appendix \ref{appendix} contains list of sequence version of functions and sequence duals covered in this paper.

\newpage 
\section{Definitions}\label{Sec:definition}
In this section, we show how to represent a sequence and definition of equivalence, sum, subtraction, multiplication, division, scalar addition and multiplication, inverse, insertion, differential, and integral for sequence.

\subsection{Sequence}
In this paper, we represent a sequence that is ordered natural numbers, or real numbers, or complex numbers, mapped from non-negative integers as follows.
  \begin{align}
    \{ a_n \} = a_0, a_1, a_2, \cdots.  
  \end{align}

\subsection{Equivalence}
First, we define the equivalence of two sequences as follows. 
\begin{definition} [Equivalence]
For any sequences $\{ a_n \}$ and $\{ b_n \}$, $\{ a_n \} = \{ b_n \}$ if and only if for any non-negative integer $n$,     
  \begin{align}
    a_n = b_n.
  \end{align}
\end{definition}

\subsection{Fundamental operations}
In this subsection, we define fundamental operations for a sequence or two sequences. First, we define shift operator $\mathcal{S}_k$ as follows.
\begin{definition} [Shift operator]
For any sequence $\{ a_n \}$ and non-negative integer $k$, 
  \begin{align}
    \mathcal{S}_k a_n = a_{n+k}, \ \ 
    \mathcal{S}_k \{ a_n \} = \{ a_{n+k} \}.
  \label{shift}
  \end{align}
\end{definition}

The sum ($+$), subtraction ($-$), multiplication ($\times$), and division ($/$) of two sequences are naturally defined as follows.
\begin{definition} [Sum, Subtraction, Multiplication, Division]
For any sequences $\{ a_n \}$ and $\{ b_n \}$, 
  \begin{align}
     \{ a_n \} + \{ b_n \} &= \{ a_n + b_n \} ,
  \\
     \{ a_n \} - \{ b_n \} &= \{ a_n - b_n \} ,
  \\
     \{ a_n \} \times \{ b_n \} &= \{ a_n \times b_n \} .
  \end{align}
For any sequence $\{ a_n \}$ and $\{ b_n \}$ that does not include $0$, 
  \begin{align}
     \{ a_n \} / \{ b_n \} = \{ a_n / b_n \} .
  \end{align}
\end{definition}

Scalar addition and scalar multiplication are also naturally defined as follows.
\begin{definition} [Scalar addition and multiplication]
For any sequence $\{ a_n \}$ and a real number $\alpha$,
  \begin{align}
     \{ a_n + \alpha \} &= \{ a_n \} + \alpha,
  \\
     \{ \alpha  a_n \} &= \alpha  \{ a_n \} .
  \end{align}
\end{definition}

Inverse of a given sequence $\{ a_n \}$ that does not include 0 is defined as follows.
\begin{definition} [Inverse]
We represent inverse of $\{ a_n \}$ that does not include 0 by $\{ a_n \}^{-1}$ and define $\{ a_n \}^{-1}$ as
  \begin{align}
     \{ a_n \}^{-1} = a_0^{-1}, \ a_1^{-1}, \cdots.
  \end{align}
\end{definition}

Also, we define insertion operator $I_{\alpha}$ for a real number $\alpha$ as follows.
\begin{definition} [Insertion]
  \begin{align}
    I_{\alpha} \{ a_n \} = \alpha, a_0, a_1, \cdots.
  \label{insertion}
  \end{align}
\end{definition}
Note that the insertion operator is relate to the shift operator as follows.
  \begin{align}
    \mathcal{S}_1 I_{\alpha} \{ a_n \} = \{ a_n \}.
  \end{align}

\subsection{Differential and Integral}
In this subsection, we define differential and integral of sequence. First, we propose two definitions of the differential of sequence: right differential and left differential. The right differential $\mathcal{D}_R \{ a_n \}_R$ of a sequence $\{ a_n \}$ is defined as follows.
\begin{definition} [Right Differential]
For any sequence $\{ a_n \}$,
  \begin{align}
    \mathcal{D}_R \{ a_n \} = \{ a_{n+1} - a_n \} .
  \label{DR}
  \end{align}
\end{definition}
We also use the left differential of $\mathcal{D}_L \{ a_n \}_L$ of a sequence $\{ a_n \}$.
\begin{definition} [Left differential]
For any sequence $\{ a_n \}$,
  \begin{align}
    \mathcal{D}_L \{ a_n \} = \{ a_n - a_{n-1} \} .
  \end{align}
\end{definition}
To make the above definition well defined, $a_{-1}$ is needed. We add an appropriate $a_{-1}$ to $\{ a_n\}$ if necessary.

We remark that the following relation of the right differential and the left differential holds.
  \begin{align} 
    \mathcal{D}_R = \mathcal{S}_1 \mathcal{D}_L.
  \end{align}

We also define the right integral $\mathcal{I}_R^{\alpha} \{ a_n \}$ of a sequence $\{ a_n \}$ is defined as follows.
\begin{definition} [Right Integral]
For any sequence $\{ a_n \}$ and an integral constant that is a real number, $\alpha$,
  \begin{align}
    \mathcal{I}_R^\alpha \{ a_n \} = \left\{ \sum_{k=0}^{n} a_k \right\} + \alpha.
  \end{align}
\end{definition}
In addition, we define the left integral $\mathcal{I}_L^\alpha \{ a_n \}$ of a sequence $\{ a_n \}$ as follows.
\begin{definition} [Left Integral]
For any sequences $\{ a_n \}$ and an integral constant that is a real number, $\alpha$,
  \begin{align}
    \mathcal{I}_L^\alpha \{ a_n \} = \left\{ \sum_{k=0}^{n-1} a_k \right\} + \alpha,
  \end{align}
where we assumed $ \sum_{k=0}^{-1} a_k = 0$.
\end{definition}

We remark that the following relation of the right differential and the left differential holds.
  \begin{align}
    \mathcal{I}_R^\alpha = \mathcal{S}_1 \mathcal{I}_L^\alpha.
  \end{align}

\begin{figure}[h]
 \centering
 \includegraphics[width=10.0cm]{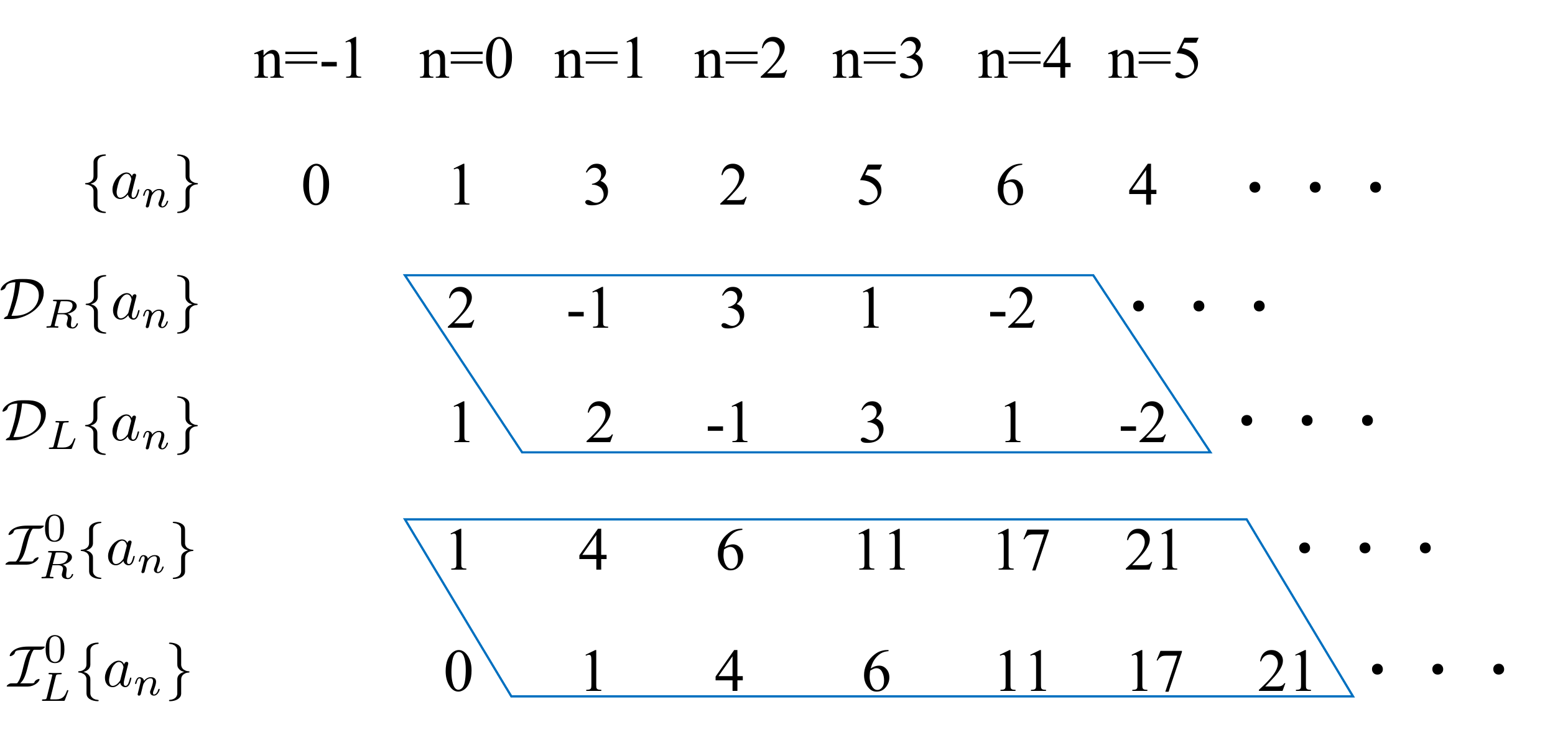}
 \caption{The left/right differential and integral of a sequence $\{a_n\}= 1,3,2,5,6,4,\cdots$ and $a_{-1}=0$.}
 \label{dif_int_ex}
\end{figure}

One can check that the right/left integral is the inverse operation of the left/right differential as follows. 
  \begin{align}
    \mathcal{D}_L  \mathcal{I}_R^{\alpha} \{ a_n \} &= \left\{ \sum_{k=0}^{n} a_k - \sum_{k=0}^{n-1} a_k \right\} = \{ a_n \},
\\
  \mathcal{I}_R^\alpha \mathcal{D}_L \{ a_n \} &= \left\{ \sum_{k=0}^{n} (a_k - a_{k-1}) \right\} +  \alpha = \{ a_n \}  - a_{-1} +  \alpha,
\\
  \mathcal{D}_R  \mathcal{I}_L^\alpha \{ a_n \} &= \left\{ \sum_{k=1}^{n} a_k - \sum_{k=1}^{n-1} a_k \right\} = \{ a_{n} \} ,
\\
  \mathcal{I}_L^\alpha \mathcal{D}_R \{ a_n \} &= \left\{ \sum_{k=0}^{n-1} (a_{k+1} - a_k) \right\} +  \alpha = \{ a_{n} - a_0 \} +  \alpha =  \{ a_n \} - a_0 + \alpha.
  \end{align}
There exist the relations among shift operator, the left/right integral, and the left/right differential as follows.
  \begin{align}
    \mathcal{D}_R  \mathcal{I}_R^\alpha \{ a_n \} &= \left\{ \sum_{k=0}^{n+1} a_k - \sum_{k=0}^{n} a_k \right\} =  \{ a_{n+1} \} = \mathcal{S}_1 \{ a_n \}, 
\\
    \mathcal{I}_R^\alpha \mathcal{D}_R \{ a_n \} &= \left\{ \sum_{k=0}^{n} (a_{k+1} - a_k) \right\} = \{ a_{n+1} - a_0 \} +  \alpha =  \mathcal{S}_1 \{ a_n \} - a_0 + \alpha, 
\\
    \mathcal{S}_1 \mathcal{D}_L  \mathcal{I}_L^\alpha \{ a_n \} &= \left\{ \sum_{k=0}^{n} a_k - \sum_{k=0}^{n-1} a_k \right\} =  \{ a_{n} \}, \\
    \mathcal{S}_1 \mathcal{I}_L^\alpha \mathcal{D}_L \{ a_n \} &= \left\{ \sum_{k=0}^n (a_k - a_{k-1}) \right\} = \{ a_n - a_0 \}  +  \alpha =   \{ a_{n} \} - a_0 + \alpha.
  \end{align}
In Fig.~\ref{dif_int_ex}, we show how the left/right differential and integral act on a sequence, $\{a_n\}=1,3,2,5,6,4,\cdots$ and $a_{-1}=0$.

\newpage  
\section{Fundamental sequence}\label{Sec:fundamental}
In this section, we define two kinds of identity sequence: additive identity sequence and multiplicative identity sequence, and sequence version of constant function. Remembering the additive operation defined in the previous section, the additive identity sequence $\{ a_n [0] \}$ should be defined as follows.
\begin{definition} [Sequence version of 0]
  \begin{align}  
    \{ a_n [0] \}  = \{ 0 \} = 0, \ 0, \ 0, \cdots.
  \end{align}
\end{definition}
Then, any sequence $\{ a_n \}$ satisfies the relations
  \begin{align}
    &\{ a_n \} + \{ a_n  [0] \}  = \{ a_n [0] \}  + \{ a_n \} = \{ a_n \},
  \\
    &\{ a_n \} - \{ a_n  [0] \}  =  \{ a_n \},
  \\
    &\{ a_n [0] \}  - \{ a_n \}  = \{ -a_n \}  = -\{ a_n \}.
  \end{align}
In addition, the multiplicative identity sequence $\{ a_n [1] \}$ should be defined as follows.
\begin{definition} [Sequence version of 1]
  \begin{align}  
    \{ a_n [1] \}  = \{ 1 \} = 1, \ 1, \ 1, \cdots.
  \end{align}
\end{definition}
Then, for any sequence $\{ a_n \}$,
  \begin{align}
    &\{ a_n \} \times \{ a_n [1] \}   = \{ a_n [1] \}  \times \{ a_n \} = \{ a_n \},
  \\
    &\{ a_n \} / \{ a_n [1] \}   = \{ a_n \},
  \end{align}
and for any sequence $\{ a_n \}$ that does not include 0,
  \begin{align}
    &\{ a_n [1] \}  / \{ a_n \} =  \{ a_n \}^{-1}.
  \end{align}
Generally, sequence version of constant function $\{ a_n [a] \}$ where $a$ is a real number is given by
\begin{definition} [Sequence version of $a$]
  \begin{align}
    \{ a_n [a] \} = \{ a \} = a, a, a, \cdots.
  \end{align}
\end{definition}
The sequence version of constant function $\{ a_n [a] \}$ satisfies the following relations for any sequence $\{ a_n \}$.
  \begin{align}
    &\{ a_n \} + \{ a_n [a] \} = \{ a_n [a] \} + \{ a_n \} = \{ a_n + a \} = \{ a_n \} + a,
  \\
    &\{ a_n \} - \{ a_n  [a] \}  =  \{ a_n - a \} = \{ a_n \} -a,
  \\
    &\{ a_n [a] \}  - \{ a_n \}  = \{ a - a_n \} = a - \{ a_n \},
  \\
    &\{ a_n \} \times \{ a_n [a] \} = \{ a_n [a] \} \times \{ a_n \} = \{ a a_n \} = a \{ a_n \},
  \\
    &\{ a_n \} / \{ a_n [a] \}   = \{ a_n / a \} =  \{ a_n \} / a,
  \label{constant_div_1}
  \\
    &\{ a_n [a] \}  / \{ a_n \} =  \{ a / a_n \} = a \{ a_n \}^{-1},
  \label{constant_div_2}
  \end{align}
where we assumed $a \neq 0$ in Eq.~(\ref{constant_div_1}) and $a_n \neq 0$ for any non-negative integer $n$ in Eq.~(\ref{constant_div_2}). In this paper, we use the following notation to represent a sequence version of a given function $f(x)$,
\begin{definition}
  \begin{align}
      \{ a_n [f(x)] \}= a_0 [f(x)], a_1 [f(x)], \cdots.
  \end{align}
\end{definition}

\newpage
\section{Power function}\label{Sec:power}
In this section, we define the sequence version of power function. First, we define $\{ a_n \}  [x]$ that is the sequence version of $x$ as follows. 
\begin{definition} [Sequence version of $x$]
  \begin{align}  
    \{ a_n  [x] \}   = \{ n \} = 0, \ 1, \ 2, \cdots.
  \end{align}
\end{definition}
Generally, we define the sequence version of power function as follows. 
\begin{definition} [Sequence version of $x^a$]
For a real number, $a$,
  \begin{align}  
    \{ a_n [x^a] \}   = \{ n^a \} = 0^a, \ 1^a, \ 2^a, \cdots.
  \end{align}
\end{definition}
Note that for real numbers $a$, $b$, the following formulae are analogous to the law of exponent.
\begin{proposition} [Sequence version of the law of exponent]
  \begin{align}  
     &\{ a_n [x^{a+b}] \}  =  \{ n^{a+b} \} = 0^{a+b}, \ 1^{a+b}, \ 2^{a+b}, \cdots 
  \nonumber\\
     &\qquad= 0^a \times 0^b, \ 1^a \times 1^b, \ 2^a \times 2^b, \cdots = \{ a_n [x^a] \} \times \{ a_n [x^b] \}  ,
  \\
     &\{ a_n [x^{a-b}] \}  =  \{ n^{a-b} \} = 0^{a-b}, \ 1^{a-b}, \ 2^{a-b}, \cdots 
  \nonumber\\
     &\qquad= \frac{0^a}{0^b}, \ \frac{1^a}{1^b}, \ \frac{2^a}{2^b}, \cdots = \{ a_n [x^a] \} / \{ a_n [x^b] \}  ,
  \\
     &\{ a_n [x^{ab}] \}  = \{ n^{ab} \} = 0^{ab}, 1^{ab}, 2^{ab}, \cdots
  \nonumber\\
     &\qquad= (0^a)^b, (1^a)^b, (2^a)^b, \cdots = \{ a_n [x^{a}] \}^b. 
  \end{align}
\end{proposition}

By using power function, one has the following theorem called Binomial theorem.
  \begin{align}
    (1 + x)^n = \sum_{k=0}^n {}_nC_k x^k.
  \end{align}

\begin{figure}[h]
 \centering
 \includegraphics[width=10.0cm]{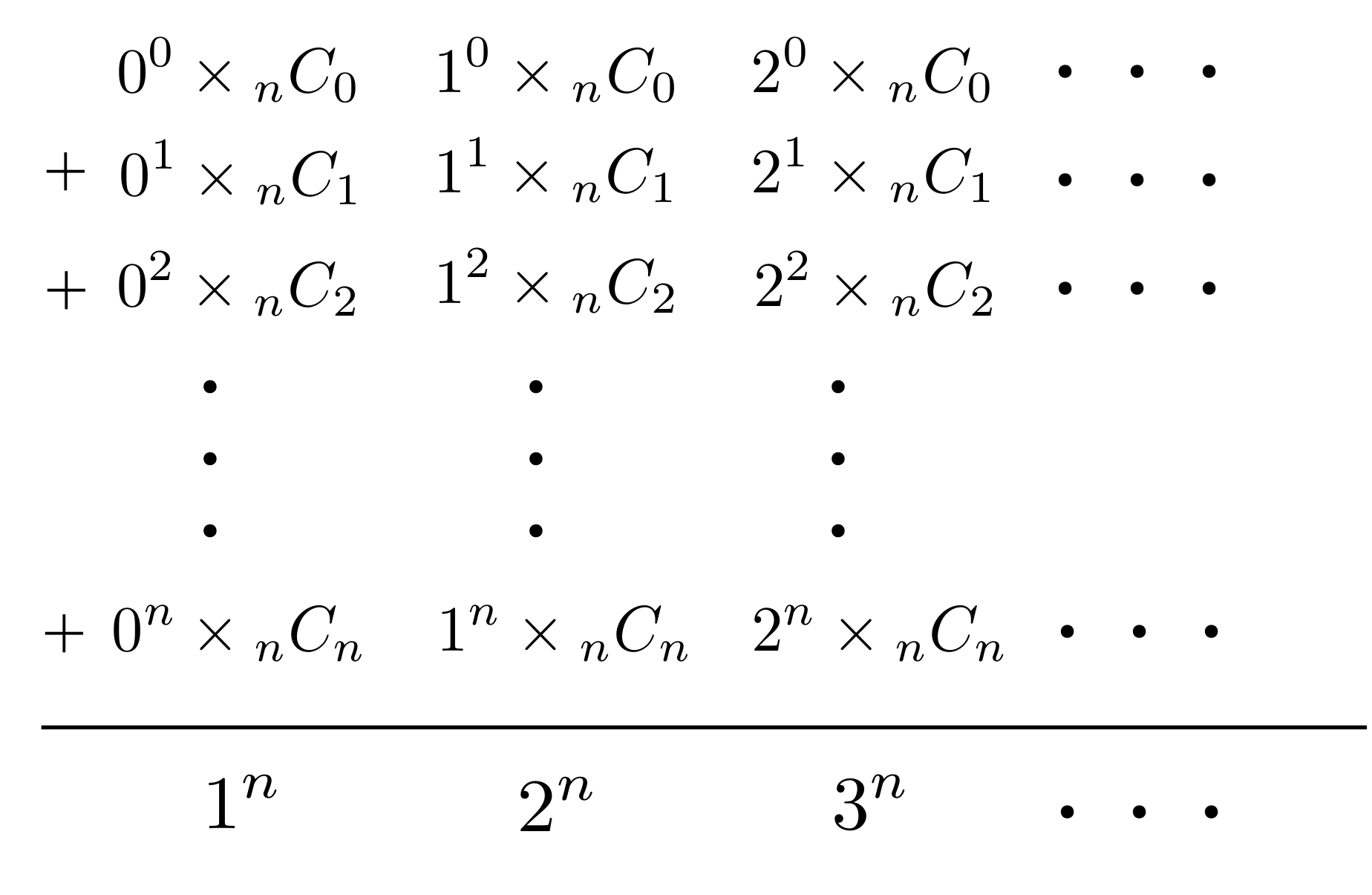}
 \caption{Diagrammatic representation of sequence version of binomial theorem.}
 \label{binomial}
\end{figure}

Then, we find the sequence version of the binomial theorem as follows (Fig.~\ref{binomial}).
\begin{theorem} [Sequence version of binomial theorem]
For a non-negative number $n$,
  \begin{align}
    &(\{ a_n [1]\} + \{ a_n [x]\})^n 
  \nonumber\\
    &\qquad= {}_nC_0 \{ a_n [1] \}^n \{ a_n [x] \}^0 + {}_nC_1 \{ a_n [1] \}^{n-1} \{ a_n [x] \}^1 + \cdots + {}_nC_n \{ a_n [1] \}^0 \{ a_n [x] \}^n
  \\
    &\qquad= \sum_{k=0}^n {}_nC_k \{ a_n [1] \}^{n-k} \{ a_n [x] \}^k.
  \end{align}
\end{theorem}

\newpage
\section{Exponential function}\label{Sec:exponential}
In this section, we define the sequence versions of exponential function by focusing on differential property of the exponential function.
\subsection{$e^x$}
\begin{figure}[h]
 \centering
 \includegraphics[width=7.5cm]{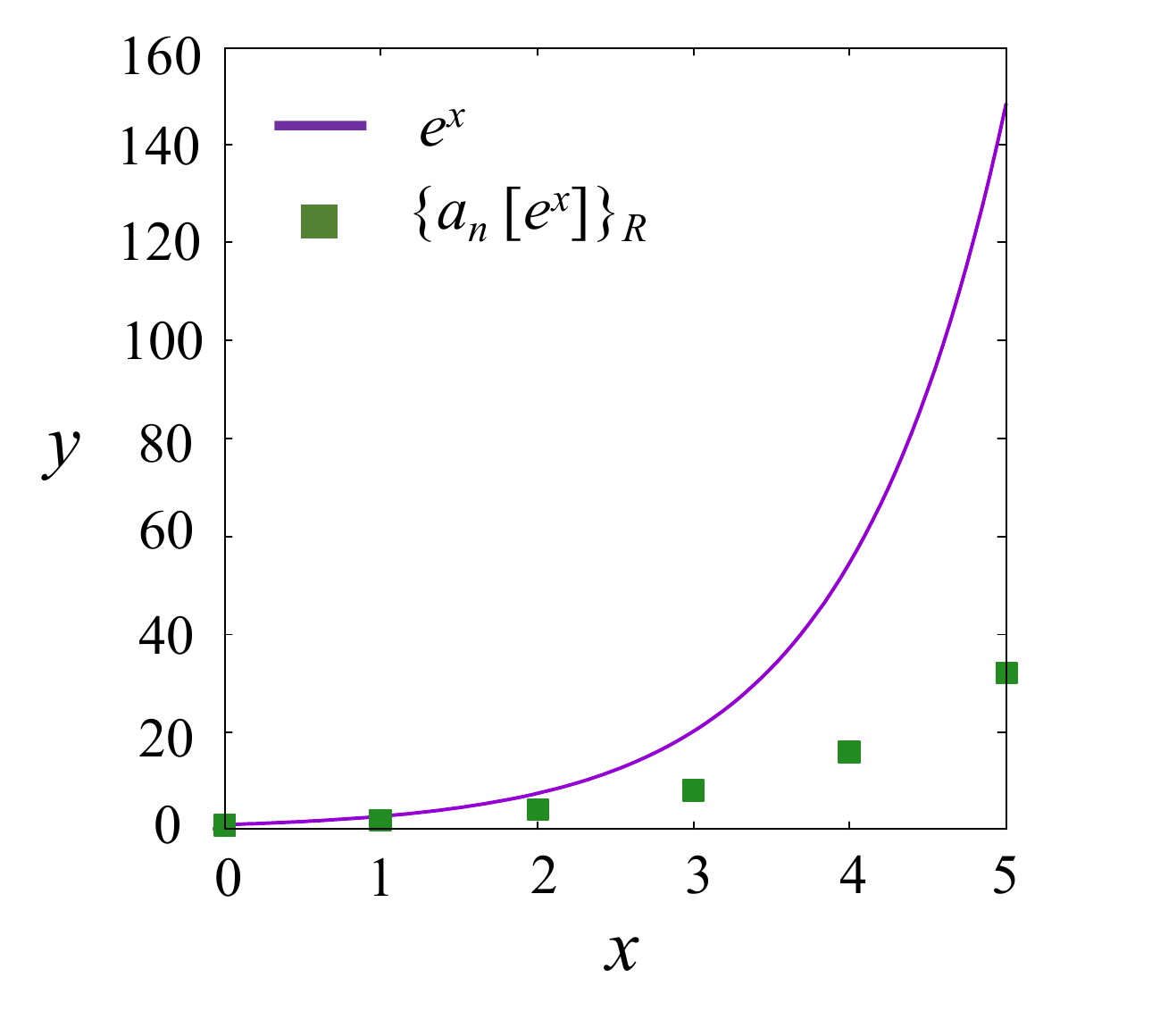}
 \caption{Difference between $e^x$ and $\{a_n [e^x]\}_R$.}
 \label{seq_ex}
\end{figure}

The exponential function satisfies the following formula for differential.
  \begin{align}
    \frac{d}{dx} e^x = e^x.
  \end{align}
Then, we find a sequence $ \{ a_n [e^x] \}_R$ that satisfies the following formula\footnote{In Remark~\ref{remark_L}, we define  $ \{ a_n [e^x] \}_L$ by using $\mathcal{D}_L$}.
  \begin{align}
    \mathcal{D}_R \{ a_n [e^x] \}_R =  \{ a_n [e^x] \}_R.
  \label{def_ex_DR}
  \end{align}  
The $\{ a_n [e^x] \}_R$ can be uniquely determined except the initial term as follows.
  \begin{align}
    &a_{n+1} - a_n  =  a_n, 
\\
    &a_{n+1} = 2a_n.
  \end{align} 
Then, we define $\{ a_n [e^x] \}_R$ as follows because $e^0=1$.
\begin{definition} [Right sequence version of exponential function]
  \begin{align}  
     \{ a_n [e^x] \}_R = \{ 2^n \} = 1, \ 2, \ 4, \cdots.
  \label{def_ex1}
  \end{align}
\end{definition}
Figure \ref{seq_ex} shows difference between the usual exponential function ($e^x$) and the right sequence version of the exponential function ($\{a_n[e^x]\}_R$).

We remark that the difference between $d/dx$ and $\mathcal{D}_R$ is manifested by the difference between the $e$ and the $2$. Then, we show sequence version of the Maclaurin series for $e^x$. By using the Maclaurin series, the $e^x$ is expressed as follows. 
  \begin{align}
    e^x &= \sum_{n=0}^{\infty} \frac{x^n}{n!} = \frac{1}{0!} 1 + \frac{1}{1!} x + \frac{1}{2!}x^2 + \frac{1}{3!}x^3 + \cdots 
\nonumber\\
    &= 1 + x  + \int_0^x x + \int_0^x \int_0^x x + \int_0^x \int_0^x \int_0^x x + \cdots .
  \end{align}
Because we have already defined the sequence version of $1$, $x$, $e^x$, and the left integral of a sequence, then, we have the following left sequence version of the Maclaurin series for $\{a_n[e^x]\}_R$ as follows.
\begin{proposition} [Left sequence version of Maclaurin series for right sequence version of exponential function 1]
  \begin{align}
  \{ a_n [e^x] \}_R  = \{ a_n  [1] \} + \{ a_n  [x] \} + \mathcal{I}_L^0 \{ a_n [x] \}  +  \mathcal{I}_L^0 \mathcal{I}_L^0 \{ a_n  [x] \} + \mathcal{I}_L^0 \mathcal{I}_L^0 \mathcal{I}_L^0 \{ a_n [x] \} + \cdots.
  \end{align}
\label{exponential}
\end{proposition}
Each term in the r.h.s.~in the above equation is given by
  \begin{align}
    \mathcal{I}_L^0 \{ a_n [x] \} &= 0, 0, 1, 3, 6, 10, \cdots,
  \\
    \mathcal{I}_L^0  \mathcal{I}_L^0 \{ a_n [x] \} &= 0, 0, 0, 1, 4, 10, 20, \cdots,
  \\
    \mathcal{I}_L^0  \mathcal{I}_L^0  \mathcal{I}_L^0 \{ a_n  [x] \} &= 0, 0, 0, 0, 1, 5, 15, 35, \cdots,
  \end{align}
etc. Then, we introduce the following definition.
\begin{definition} [Left sequence version of $x^n/n!$]
For an integer $k \geq2$,
  \begin{align}
    \{ a_n [x^k/k!] \}_L =  \underbrace{\mathcal{I}_L^0 \cdots \mathcal{I}_L^0}_{\rm k-1 \ times} \{ a_n [x] \}.
  \end{align}
\end{definition} 
We remark $ \{ a_n [x^2] \}  \neq \{ a_n [x^2/2] \}_L  + \{ a_n [x^2/2] \}_L $. Instead, the following formula holds.
  \begin{align}
    \{ a_n [x^2] \} =  \{ a_n  [x^2/2] \}_L + \mathcal{S}_{1} \{ a_n [x^2/2]  \}_L.
  \end{align}
Generally, we obtain the following proposition.
\begin{proposition} [Relation between left sequence version of $x^n/n!$ and sequence version of $x^n$] 
For a natural number $n$,
  \begin{align} 
    \{ a_n [x^n] \} =  \sum_{k=0}^{n-1} A(n, k) \mathcal{S}_k \{ a_n [x^n/n!] \}_L,
  \label{relation_Euler}
  \end{align}
where $A(n, k)$ is the Eulerian number satisfying the following formula for $n=1,2,3,\cdots$ and $k=0,1,2,\cdots,n$ \cite{Euler:number, Euler:number2, Euler:number3}.
  \begin{align}
    A(n, k) = \sum_{l=0}^{k+1} (-1)^l 
      \begin{pmatrix}
      n+1
      \\
      l
      \end{pmatrix}
    (k+1-l)^n.
  \end{align}
\label{Eulerian_number}
\end{proposition}
The Eulerian number represents the number of permutations for $n$ elements with $k$ raises. For example, $A(3,2)=1$ because there exist only one permutation for 3 elements with 2 raises: $(1,2,3)$. Note that $3-2=+1$ and $2-1=+1$.

Then, we have the following formula.
\begin{proposition} [Left sequence version of Maclaurin series for right sequence version of exponential function 2] 
  \begin{align}
     \{ a_n [e^x] \}_R = \sum_{n=0}^{\infty} \{ a_n [x^n/n!] \}_L.
  \end{align}
\label{exponential2}
\end{proposition}

By using the right integral of a sequence, we have the following right sequence version of the Maclaurin series for $\{a_n[e^x]\}_R$ as follows.
\begin{proposition} [Right sequence version of Maclaurin series for right sequence version of exponential function 1]
  \begin{align}
      \{ a_n [e^x] \}_R  &= \{ a_n  [1] \} + \frac{1}{2} \{ a_n  [x] \} + \frac{1}{4} \mathcal{I}_R^0 \{ a_n [x] \} 
  \nonumber\\
      &\qquad +  \frac{1}{8}  \mathcal{I}_R^0 \mathcal{I}_R^0 \{ a_n  [x] \} + \frac{1}{16}  \mathcal{I}_R^0 \mathcal{I}_R^0 \mathcal{I}_R^0 \{ a_n [x] \} + \cdots.
  \end{align}
\label{exponential3}
\end{proposition}
Also, we have
\begin{proposition} [Right sequence version of Maclaurin series for right sequencne version of exponential function 2] 
  \begin{align}
     \{ a_n [e^x] \}_R = \sum_{n=0}^{\infty} \frac{1}{2^n} \{ a_n [x^n/n!] \}_R.
  \end{align}
\label{exponential4}
\end{proposition}
Here, we defined $\{ a_n [x^n/n!] \}_R$ as follows.
\begin{definition} [Right sequence version of $x^n/n!$]
For an integer $k \geq2$,
  \begin{align}
    \{ a_n [x^k/k!] \}_R =  \underbrace{\mathcal{I}_R^0 \cdots \mathcal{I}_R^0}_{\rm k-1 \ times} \{ a_n [x] \} = \mathcal{S}_{k-1} \{ a_n [x^k/k!] \}_L .
  \end{align}
\end{definition} 
For example, we have
  \begin{align}
    \mathcal{I}_R^0 \{ a_n [x] \} &= 0, 1, 3, 6, 10, \cdots,
  \\
    \mathcal{I}_R^0  \mathcal{I}_R^0 \{ a_n [x] \} &= 0, 1, 4, 10, 20, \cdots,
  \\
    \mathcal{I}_R^0  \mathcal{I}_R^0  \mathcal{I}_R^0 \{ a_n  [x] \} &= 0,  1, 5, 15, 35, \cdots.
  \end{align}
Also, the first three terms of the right sequence version of the Maclaurin series for $\{a_n[e^x]\}_R$ are represented as follows.
  \begin{align}
    & 1 + 0 + 0 + \cdots = 1,
  \\
    & 1 + \frac{1}{2} + \frac{1}{4} + \cdots = \sum_{n=0}^{\infty} \frac{1}{2^n} = 2,
  \\
    & 1 + \frac{2}{2} + \frac{3}{4} + \cdots = \sum_{n=0}^{\infty} \frac{n+1}{2^n} =4.
  \end{align}
In addition, we obtain the following relation between the right sequence version of $x^n/n!$ and sequence version of $x^n$.
  \begin{proposition} [Relation between right sequence version of $x^n/n!$ and sequence version of $x^n$] 
  For a natural number $n$,
    \begin{align} 
      \{ a_n [x^n] \} =  \sum_{k=0}^{n-1} A(n, k) I_0^{n-1-k} \{ a_n [x^n/n!] \}_R,
    \end{align}
  \label{relation_Euler_2}
  \end{proposition}
Note that the factors $1/2$, $1/2^2$ etc.\,are needed to construct the right sequence version of Maclaurin series for $\{ a_n  [e^x] \}_R$, and replacement of the $\int_0^x$ by the $\mathcal{I}_L^0$ does not directly lead to the right sequence version of Maclaurin series for $\{ a_n  [e^x] \}_R$. In this paper, the factors $1/2$, $1/2^2$ etc.\,frequently play the role of bridge between the world of sequence and the world of analytics.

We remark the following three points.

\begin{remark}
We show two kinds of diagrams corresponding to Prop.~\ref{exponential} in Figs.~\ref{ex_Pascal} and ~\ref{block_ex}.
\begin{figure}[h]
 \centering
 \includegraphics[width=14.5cm]{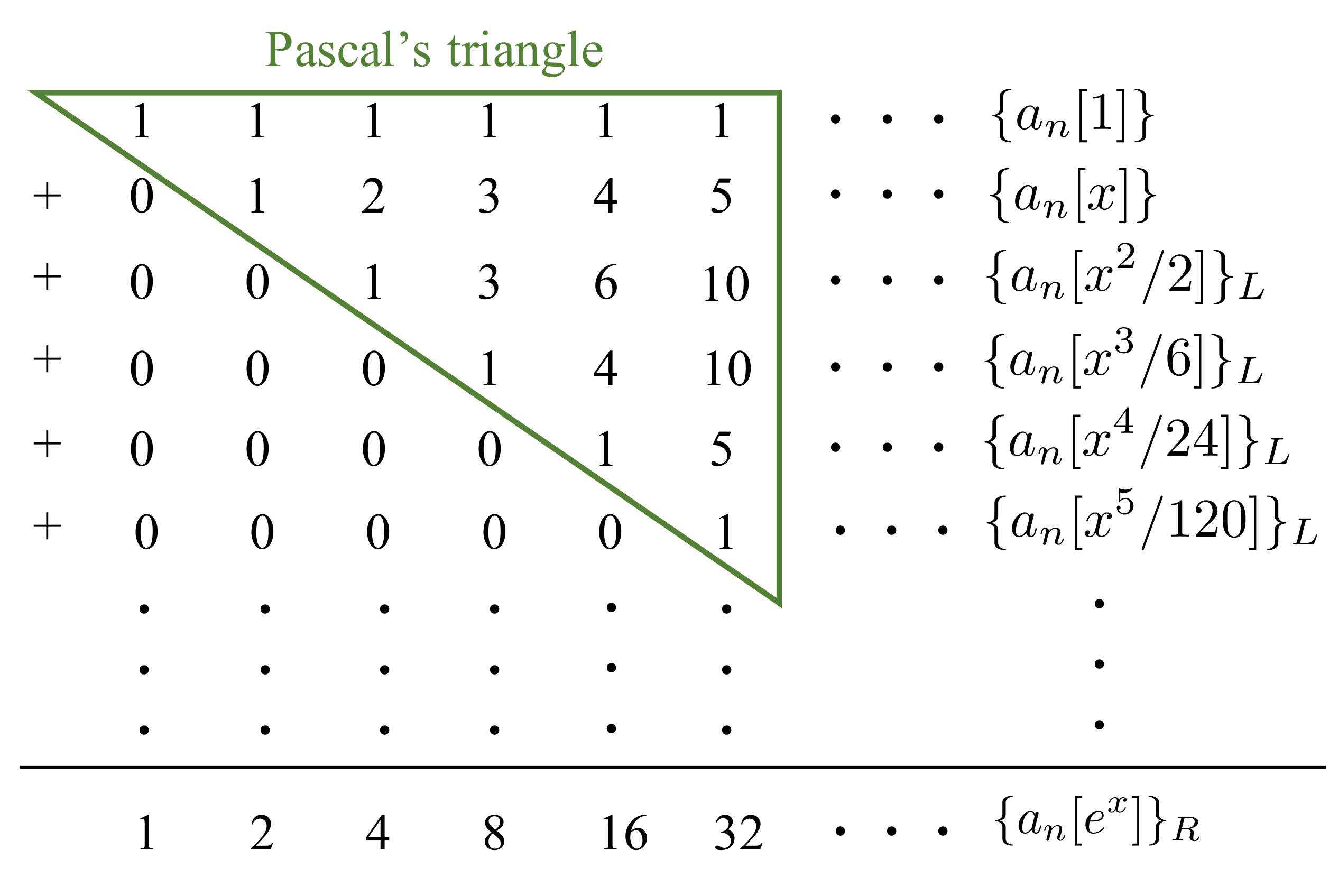}
 \caption{Diagrammatic representation of left sequence version of Maclaurin series for $\{ a_n [e^x] \}_R$ that is equivalent to Pascal's triangle.}
 \label{ex_Pascal}
\end{figure}

\begin{figure}[h]
 \centering
 \includegraphics[width=14.5cm]{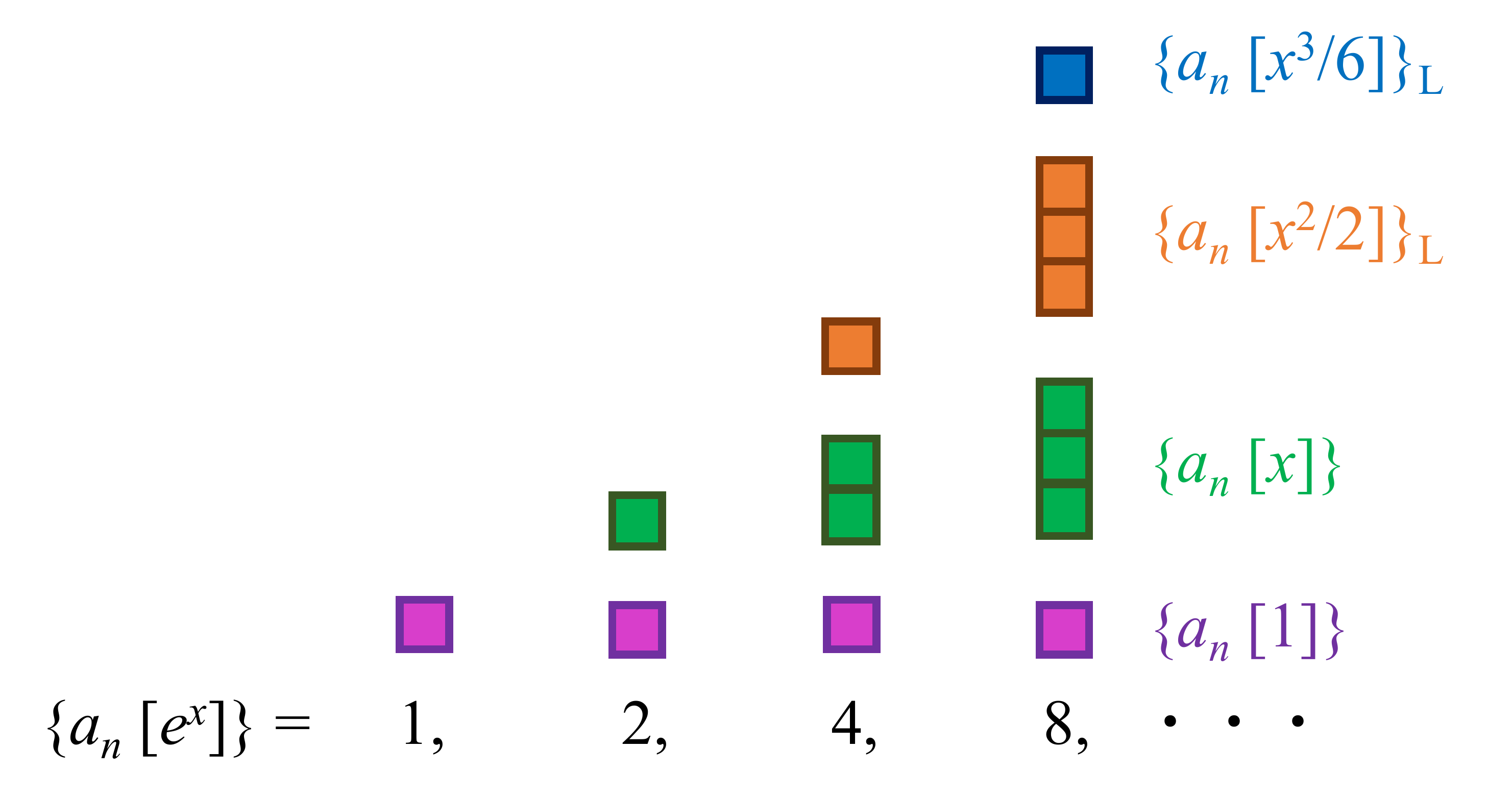}
 \caption{Diagrammatic representation of left sequence version of Maclaurin series for $\{ a_n [e^x] \}_R$ by using blocks.}
 \label{block_ex}
\end{figure}
\end{remark}

\begin{remark}\label{remark_L}
We used $\mathcal{D}_R$ to define $\{ a_n[e^x] \}_R$ in Eq.~(\ref{def_ex_DR}). The $\{ a_n[e^x] \}_L$ should satisfy the following properties for any non-negative integer $n$.
  \begin{align}
    a_n - a_{n-1} &= a_n,
  \\
    a_{n-1} & = 0.
  \end{align}
Then, 
  \begin{align}
    \{ a_n[e^x] \}_L = 0, 0, \cdots.
  \end{align}
\end{remark}

\begin{remark}
In Fig.~\ref{Eulerian_3}, Prop.~\ref{Eulerian_number} for $n=3$ is represented.
\begin{figure}[h]
 \centering
 \includegraphics[width=14.5cm]{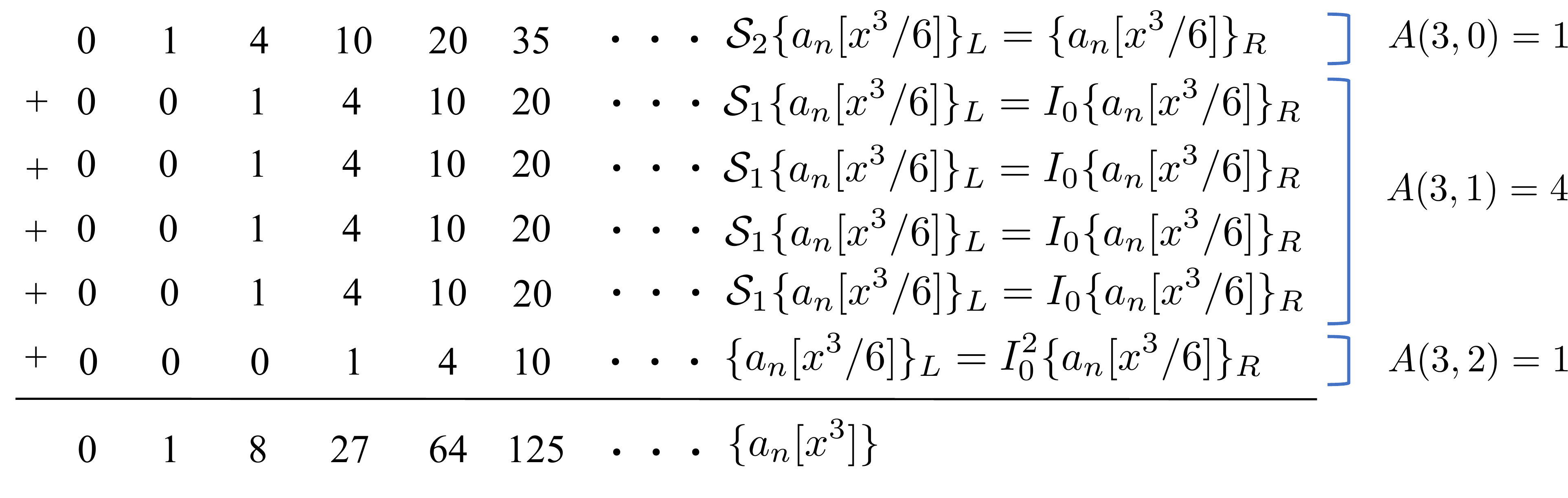}
 \caption{Diagrammatic representation of the relation among $\{ a_n [x^n] \}$, $\{ a_n [x^n/n!] \}_L$, and $\{ a_n [x^n/n!] \}_R$ with the Eulerian number $A(n,k)$ for $n=3$.}
 \label{Eulerian_3}
\end{figure}
\end{remark}

\subsection{$e^{-x}$}
Next, we define the sequence version of $e^{-x}$. The $e^{-x}$ satisfies the following formula.
  \begin{align}
    \frac{d}{dx} e^{-x} = -e^{-x}.
  \end{align}
Then, we find a sequence $ \{ a_n [e^{-x}] \}_L$ that satisfies the following formula.
  \begin{align}
    \mathcal{D}_L \{ a_n [e^x] \}_L = - \{ a_n [e^x] \}_L.
  \end{align} 
The $\{ a_n [e^{-x}]\}_L$ can uniquely be determined except the initial term as follows.
  \begin{align}
    &a_{n} - a_{n-1}  =  -a_n, 
\\
    &2a_{n} = a_{n-1}.
  \end{align} 
Then, we define $\{ a_n [e^{-x}] \}_L = a_0 [e^{-x}]_L, a_1 [e^{-x}]_L, \cdots$ as follows because $e^{-0}=1$.
\begin{definition} [Left sequence version of negative exponential function]
  \begin{align}  
     \{ a_n [e^{-x}] \}_L = \{ 2^{-n} \} = 1, \ \frac{1}{2}, \ \frac{1}{4}, \cdots,
  \label{def_ex2}
  \end{align}
\end{definition}
where we assumed $ a_{-1} [e^{-x}]_L =2$.

By combining the right sequence version of exponential function and the left sequence version of negative exponential function, we have the following inverse relation.
\begin{proposition} [Inverse relation between right sequence version of exponential function and left sequence version of negative exponential function]
  \begin{align}
    \{ a_n [(e^x)] \}_R^{-1} = \{ a_n [e^{-x}] \}_L.
  \label{inverse_ex}
  \end{align}
\end{proposition}
This is analogous to the following relation.
  \begin{align}
    (e^x)^{-1} = e^{-x} .
  \end{align}

As with the left sequence version of Maclaurin series for $\{ a_n  [e^{x}] \}_R$, the left sequence version of Maclaurin series for $\{ a_n  [e^{-x}] \}_L$ can be obtained. First, the Maclaurin series for the $e^{-x}$ is expressed as follows.
  \begin{align}
    e^{-x} &= \sum_{n=0}^{\infty} \frac{(-1)^n x^n}{n!} = \frac{1}{0!} 1 - \frac{1}{1!} x + \frac{1}{2!}x^2 - \frac{1}{3!}x^3 + \cdots 
\\
    &= 1 - x  + \int_0^x x - \int_0^x \int_0^x x + \int_0^x \int_0^x \int_0^x x - \cdots .
  \end{align}
Then, the left sequence version of Maclaurin series for $\{ a_n  [e^{-x}] \}_L$ is obtained as follows.
\begin{proposition} [Left sequence version of Maclaurin series for left sequence version of negative exponential function 1]
  \begin{align}
     \{ a_n [e^{-x}] \}_L &=  \{ a_n [1] \}  - \frac{1}{2} \{ a_n  [x] \} + \frac{1}{2^2}  \mathcal{I}_L^0 \{ a_n [x] \}  
  \nonumber\\
     &\qquad - \frac{1}{2^3}  \mathcal{I}_L^0 \mathcal{I}_L^0 \{ a_n  [x] \} + \frac{1}{2^4} \mathcal{I}_L^0 \mathcal{I}_L^0 \mathcal{I}_L^0 \{ a_n [x] \} - \cdots.
  \end{align}
\label{exponential1-}
\end{proposition}

By using $\{ a_n [x^n/n!] \}_L$, we can rewrite Prop.~\ref{exponential1-} as follows.
\begin{proposition} [Left sequence version of Maclaurin series for left sequence version of negative exponential function 2]
  \begin{align}
    \{ a_n [e^{-x}] \}_L =  \sum_{n=0}^{\infty} \frac{(-1)^n}{2^n} \{ a_n [x^n/n!] \}_L.
  \end{align}
\label{exponential2-}
\end{proposition}

Also, by using the left integral, we obtain the following proposition of the right sequence version of Maclaurin series for $\{ a_n  [e^{-x}] \}_L$.
\begin{proposition} [Right sequence version of Maclaurin series for left sequence version of negative exponential function 1]
  \begin{align}
     \{ a_n [e^{-x}] \}_L &=  \{ a_n [1] \}  -  \{ a_n  [x] \} + \frac{1}{2^2}  \mathcal{I}_L^0 \{ a_n [x] \}  -   \mathcal{I}_L^0 \mathcal{I}_L^0 \{ a_n  [x] \}
  \nonumber\\
     &\qquad +  \mathcal{I}_L^0 \mathcal{I}_L^0 \mathcal{I}_L^0 \{ a_n [x] \} - \cdots.
  \end{align}
\label{exponential3-}
\end{proposition}

By using $\{ a_n [x^n/n!] \}_R$, we can rewrite Prop.~\ref{exponential3-} as follows.
\begin{proposition} [Right sequence version of Maclaurin series for left sequence version of negative exponential function 2]
  \begin{align}
    \{ a_n [e^{-x}] \}_L =  \sum_{n=0}^{\infty} (-1)^n \{ a_n [x^n/n!] \}_R.
  \end{align}
\label{exponential4-}
\end{proposition}

We remark that the sums in Props.~\ref{exponential3-} and \ref{exponential4-} should not be interpreted by the usual sum but the Abel sum \cite{divergent}. For example, the first term and second term of the right sequence version of Maclaurin series for $\{ a_n [e^{-x}] \}_L$ are given as follows.
  \begin{align}
    &1 - 1 + 1 - 1 + 1 - 1 + \cdots = \frac{1}{2} =  a_1 [x^n/n!]_R,
  \\
    &1 - 2 + 3 - 4 + 5 - 6 \cdots = \frac{1}{4}  = a_2 [x^n/n!]_R.
  \end{align}
The first one is called Grandi's series \cite{Grandi}, and would be one of the most famous divergent series. The right sequence version of Maclaurin series for all the terms of $\{ a_n [e^{-x}] \}_L$ are divergent except the zeroth term if the usual sum is used. However, by using the Abel sum, value of real numbers can be assigned to those divergent series. Note that the $n$th term of $\{ a_n [e^{-x}] \}_L$ can be represented as follows.
  \begin{align}
    \left. a_n [e^{-x}]_L = \frac{1}{(1+x)^n} \right|_{x=1}.
  \end{align}
If $|x| < 1$, the r.h.s.~of the above equation is given by 
  \begin{align}
    \frac{1}{(1+x)^n} &= 1 - nx + \frac{n(n+1)}{2} x^2 -  \frac{n(n+1)(n+2)}{6} x^3 + \cdots
  \nonumber\\
    &= 1 + \sum_{k=1}^{\infty} \frac{n(n+1)\cdots(n+k-1)}{k!} (-x)^k.
 \end{align}
If the above sum is interpreted by the Abel sum, the above formula holds even if $|x| = 1$.

\begin{remark}
Figure \ref{ex-_Pascal} (\ref{ex-_Pascal_2}) shows the diagram representing the right sequence version of Maclaurin series for $\{ a_n  [e^{-x}] \}_L$ (left sequence version of Maclaurin series for $\{ a_n  [e^{-x}] \}_L$). 
\begin{figure}[h]
 \centering
 \includegraphics[width=14.5cm]{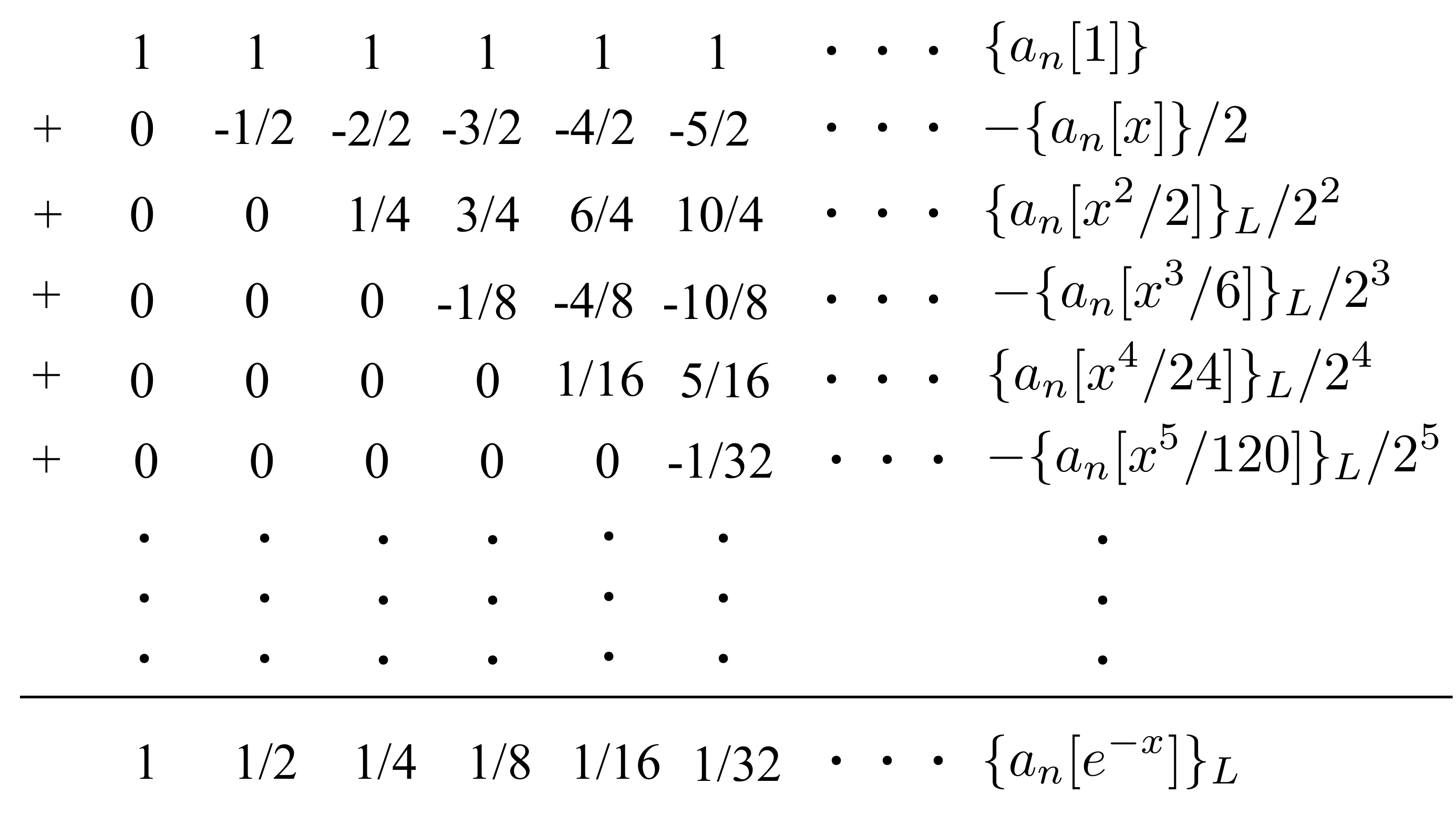}
 \caption{Diagrammatic representation of right sequence version of Maclaurin series for left sequence version of negative exponential function}
 \label{ex-_Pascal}
\end{figure}

\begin{figure}[h]
 \centering
 \includegraphics[width=14.5cm]{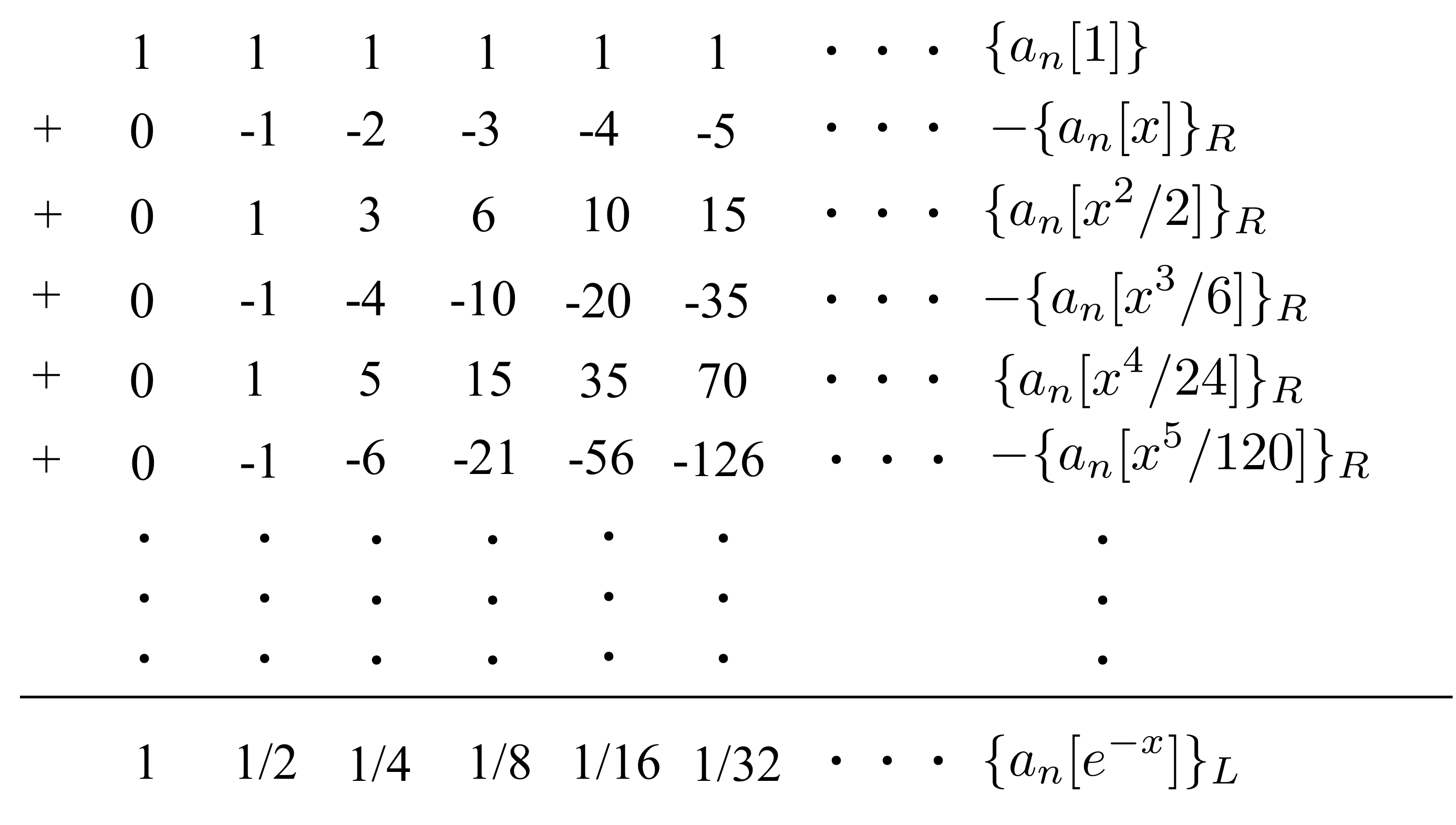}
 \caption{Diagrammatic representation of left sequence version of Maclaurin series for left sequence version of negative exponential function}
 \label{ex-_Pascal_2}
\end{figure}

\end{remark}

\begin{remark}
The $\{ a_n[e^{-x}] \}_R$ should satisfy the following properties for any non-negative integer $n$.
  \begin{align}
    a_{n+1} - a_{n} &= -a_n,
  \\
    a_{n+1} & = 0.
  \end{align}
Then, 
  \begin{align}
    \{ a_n[e^{-x}] \}_R = a, 0, 0, 0, \cdots,
  \end{align} 
where $a$ is a real number.
\end{remark}

\begin{remark}
The sequence $\{ a_n[e^{-x}] \}^{\rm nat}$ naturally corresponding to $e^{-x}$ from the viewpoint of Maclaurin series is given by
  \begin{align}
    &\{ a_n[e^{-x}] \}^{\rm nat}  =  \{ a_n [1]  \} -  \{ a_n  [x] \} +   \mathcal{I}_L^0 \{ a_n [x] \}  - \mathcal{I}_L^0 \mathcal{I}_L^0 \{ a_n  [x] \}
  \nonumber\\
     &\qquad + \mathcal{I}_L^0 \mathcal{I}_L^0 \mathcal{I}_L^0 \{ a_n [x] \} - \cdots = \{ (1-1)^n \} = 1, 0, 0, \cdots,
  \label{ex-_natural}
  \end{align}
where we assumed $0^0=1$.
\end{remark}

\subsection{$e^{\alpha x}$}
In this subsection, we generalize the result obtained in the previous two subsections, and define sequence version of $e^{\alpha x}$ for any real number $\alpha$. 
\subsubsection{$\{ a_n [e^{\alpha x}] \}_R, \ \alpha\neq-1$}
A sequence satisfying the following relation
  \begin{align}
    \mathcal{D}_R \{ a_n [e^{\alpha x}] \}_R = \alpha \{ a_n [e^{\alpha x}] \}_R.
  \end{align}  
is given by
\begin{definition} [Right sequence version of general exponential function]
  \begin{align}
    \{ a_n [e^{\alpha x}] \}_R = 1, (1+\alpha), (1+\alpha)^2, \cdots.
  \label{def_ex_gene1}
  \end{align}
\end{definition}

Then, the exponential function satisfies the law of exponent as follows for real numbers $a$ and $b$.
  \begin{align}
    e^{(a+b)x} = e^{ax} e^{bx}.
  \label{exponent_exponential_1}
  \end{align}
However, the law of exponent for $\{a_n [e^{\alpha x}] \}_R$ itself does not hold and instead, the following relation holds.
\begin{proposition} [Law of exponent for right sequence version of general exponential function]
For real numbers $\alpha$ and $\beta$, 
  \begin{align}
    &\{a_n [e^{( \alpha  + \beta  + \alpha \beta ) x}] \}_R = \{a_n [e^{ \alpha x}] \}_R \{a_n [e^{ \beta x}] \}_R .
  \end{align}
\end{proposition}

Also, the following four propositions of the sequence versions of Maclaurin series for general exponential function hold.
\begin{proposition} [Left sequence version of Maclaurin series for right sequence version of general exponential function 1]
The sequence version of Maclaurin series for $\{ a_n [e^{\alpha x}] \}_R$ is given by
  \begin{align}
    \{ a_n [e^{\alpha x}] \}_R &= \{ a_n [1] \}  + \alpha \{ a_n  [x] \} + \alpha^2 \mathcal{I}_L^0 \{ a_n [x] \}  + \alpha^3 \mathcal{I}_L^0 \mathcal{I}_L^0 \{ a_n  [x] \}
  \nonumber\\
    &\qquad + \alpha^4 \mathcal{I}_L^0 \mathcal{I}_L^0 \mathcal{I}_L^0 \{ a_n  [x] \}+ \cdots.
  \end{align}
\end{proposition}

\begin{proposition} [Left sequence version of Maclaurin series for right sequence version of general exponential function 2]
The sequence version of Maclaurin series for $\{ a_n [e^{\alpha x}] \}_R$ is given by
  \begin{align}
    \{ a_n [e^{\alpha x}] \}_R &= \sum_{n=0}^{\infty} \alpha^n \{ a_n [x^n/n!] \}_L.
  \end{align}
\end{proposition}

\begin{proposition} [Right sequence version of Maclaurin series for right sequence version of general exponential function 1]
The sequence version of Maclaurin series for $\{ a_n [e^{\alpha x}] \}_R$ is given by
  \begin{align}
    \{ a_n [e^{\alpha x}] \}_R &= \{ a_n [1] \}  + \frac{\alpha}{1 + \alpha} \{ a_n  [x] \} + \left(\frac{\alpha}{1 + \alpha}\right)^2 \mathcal{I}_R^0 \{ a_n [x] \}  +\left(\frac{\alpha}{1 + \alpha}\right)^3 \mathcal{I}_R^0 \mathcal{I}_R^0 \{ a_n  [x] \}
  \nonumber\\
    &\qquad + \left(\frac{\alpha}{1 + \alpha}\right)^4 \mathcal{I}_R^0 \mathcal{I}_R^0 \mathcal{I}_R^0 \{ a_n  [x] \}+ \cdots.
  \end{align}
\end{proposition}

\begin{proposition} [Right sequence version of Maclaurin series for right sequence version of general exponential function 2]
The sequence version of Maclaurin series for $\{ a_n [e^{\alpha x}] \}_R$ is given by
  \begin{align}
    \{ a_n [e^{\alpha x}] \}_R &= \sum_{n=0}^{\infty} \left(\frac{\alpha}{1 + \alpha}\right)^n \{ a_n [x^n/n!] \}_R.
  \end{align}
\end{proposition}

The right sequence version of Maclaurin series should be interpreted through Borel summation \cite{divergent}. The following formula usually holds if $|\alpha|<1$.
  \begin{align}
    1 + \frac{\alpha}{1+\alpha} + \left( \frac{\alpha}{1+\alpha} \right)^2 + \cdots = 1 + \alpha.
  \end{align}
However, if the r.h.s.~is interpreted by the Borel summation, for all $\alpha\neq-1$, the divergent series appeared in the right sequence version of Maclaurin series can have value of real numbers. For example, the sequence version of Maclaurin series for the $ a_1 [e^{-\frac{2}{3} x}]_R $ is given by 
  \begin{align}
    1 - 2 + 4 - 8 + 16 - 32 + \cdots = \frac{1}{3} = a_1 [e^{-\frac{2}{3} x}]_R.
  \end{align}

\subsubsection{$\{ a_n [e^{\alpha x}] \}_L, \ \alpha\neq1$}
A sequence satisfying the following relation
  \begin{align}
    \mathcal{D}_L \{ a_n [e^{\alpha x}] \}_L = \alpha \{ a_n [e^{\alpha x}] \}_L,
  \end{align}
is given by 
\begin{definition}[Left sequence version of general exponential function]
  \begin{align}
    \{ a_n [e^{\alpha x}] \}_L = 1, (1  - \alpha)^{-1}, (1 -\alpha)^{-2}, \cdots.
  \label{def_ex_gene2}
  \end{align}
\end{definition}

Then, the sequence version of the law of exponent of exponential function (Eq.~(\ref{exponent_exponential_1})) is given as follows.
\begin{proposition} [Law of exponent for left sequence version of general exponential function]
For real numbers $\alpha$ and $\beta$, 
  \begin{align}
    &\{a_n [e^{( \alpha + \beta - \alpha \beta ) x}] \}_L = \{a_n [e^{ \alpha x}] \}_L \{a_n [e^{ \beta x}] \}_L .
  \end{align}
\end{proposition}

Also, we have a generalization of the inverse relation (Eq.~\ref{inverse_ex}) as follows. 
\begin{proposition} [Inverse relation for sequence versions of general exponential function]
For a real number $\alpha$ that is not equals to $-1$,
  \begin{align}
    \{ a_n [e^{ \alpha x}] \}_R^{-1} = \{a_n [e^{ -\alpha x}] \}_L.
  \end{align}
\end{proposition}

Furthermore, the following four propositions of the sequence versions of Maclaurin series for general exponential function hold.
\begin{proposition} [Left sequence version of Maclaurin series for left sequence version of general exponential function 1]
The sequence version of Maclaurin series for $\{ a_n [e^{\alpha x}] \}_L$ is given by
  \begin{align}
    \{ a_n [e^{\alpha x}] \}_L &= \{ a_n  [1] \} + \frac{\alpha}{1-\alpha} \{ a_n  [x] \} + \left(\frac{\alpha}{1-\alpha} \right)^2 \mathcal{I}_L^0 \{ a_n [x] \}  
  \nonumber\\
    &\qquad + \left(\frac{\alpha}{1-\alpha} \right)^3 \mathcal{I}_L^0 \mathcal{I}_L^0 \{ a_n  [x] \} + \left(\frac{\alpha+1}{\alpha} \right)^4 \mathcal{I}_L^0 \mathcal{I}_L^0 \mathcal{I}_L^0 \{ a_n  [x] \}+ \cdots.
  \end{align}
\end{proposition}

\begin{proposition} [Left sequence version of Maclaurin series for left sequence version of general exponential function 2]
The sequence version of Maclaurin series for $\{ a_n [e^{\alpha x}] \}_R$ is given by
  \begin{align}
    \{ a_n [e^{\alpha x}] \}_L &= \sum_{n=0}^{\infty} \left( \frac{\alpha}{1-\alpha} \right)^n \{ a_n [x^n/n!] \}_L.
  \end{align}
\end{proposition}

\begin{proposition} [Right sequence version of Maclaurin series for left sequence version of general exponential function 1]
The sequence version of Maclaurin series for $\{ a_n [e^{\alpha x}] \}_R$ is given by
  \begin{align}
    \{ a_n [e^{\alpha x}] \}_L &= \{ a_n [1] \}  + \alpha \{ a_n  [x] \} + \alpha^2 \mathcal{I}_R^0 \{ a_n [x] \}  + \alpha^3 \mathcal{I}_R^0 \mathcal{I}_R^0 \{ a_n  [x] \}
  \nonumber\\
    &\qquad + \alpha^4 \mathcal{I}_R^0 \mathcal{I}_R^0 \mathcal{I}_R^0 \{ a_n  [x] \}+ \cdots.
  \end{align}
\end{proposition}

\begin{proposition} [Right sequence version of Maclaurin series for left sequence version of general exponential function 2]
The sequence version of Maclaurin series for $\{ a_n [e^{\alpha x}] \}_R$ is given by
  \begin{align}
    \{ a_n [e^{\alpha x}] \}_L &= \sum_{n=0}^{\infty} \alpha^n \{ a_n [x^n/n!] \}_R.
  \end{align}
\end{proposition}

\newpage
\section{Hyperbolic function}\label{Sec:hyperbolic}
\begin{figure}[h]
 \centering
 \includegraphics[width=14.5cm]{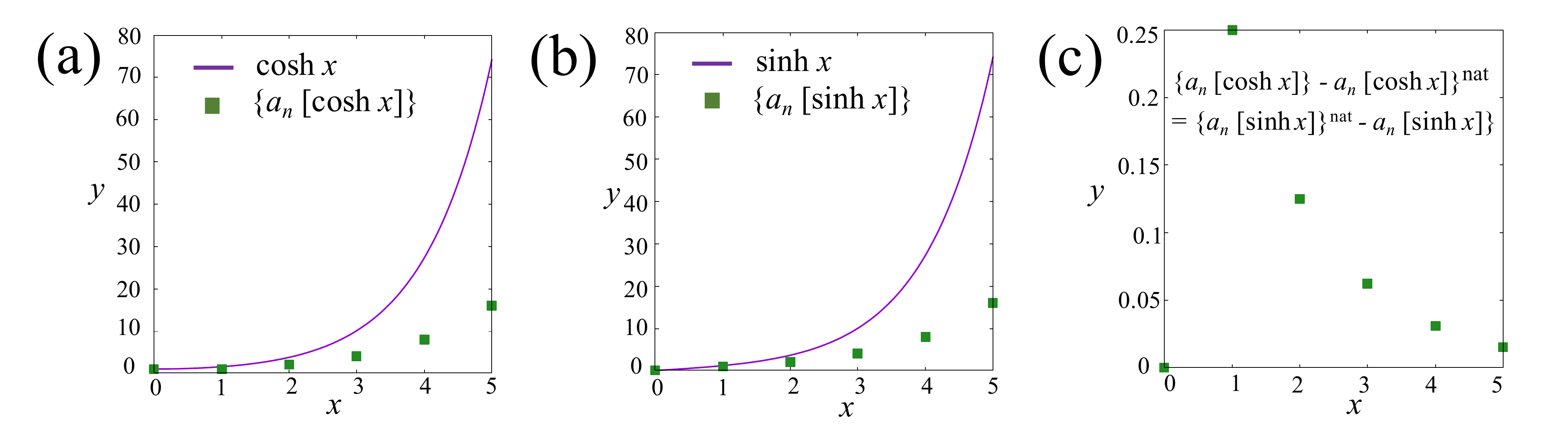}
 \caption{Difference between (a) $\cosh x$ and $\{a_n[\cosh x]\}$, (b) $\sinh x$ and $\{a_n[\sinh x]\}$, (c) $\{a_n[\cosh x]\}$ ($\{a_n[\sinh x]\}^{\rm nat}$) and $\{a_n[\cosh x]\}^{\rm nat}$ ($\{a_n[\sinh x]\}$).}
 \label{seq_coshx_sinhx}
\end{figure}

In this subsection, we define the sequence version of hyperbolic function in two ways. Because the hyperbolic cosine function and hyperbolic sine function are defined as follows.
  \begin{align}
    \cosh x = \frac{e^x + e^{-x}}{2} = \sum_{n=0}^{\infty} \frac{x^{2n}}{(2n)!},
  \\
    \sinh x = \frac{e^x - e^{-x}}{2} = \sum_{n=0}^{\infty} \frac{x^{2n+1}}{(2n+1)!},
  \end{align}
and we have already defined the sequence version of the $e^x$ in Eq.~(\ref{def_ex1}) and $e^{-x}$ in Eq.~(\ref{def_ex2}), then, we can obtain sequence version of the hyperbolic cosine function and hyperbolic sine function as follows.
\begin{definition} [Sequence version of hyperbolic function]
  \begin{align}
    \{ a_n [\cosh x] \} = \frac{1}{2} [\{ a_n [[e^{x}] ] \}_R + \{ a_n [e^{-x}] \}_L ] = \frac{1}{2} \{ 2^{n} + 2^{-n} \} = 1, \frac{5}{4}, \frac{17}{8}, \frac{65}{16}, \cdots,
  \label{cosh_seq}
  \\
    \{ a_n [\sinh x] \} = \frac{1}{2} [\{ a_n [[e^{x}] ] \}_R - \{ a_n [e^{-x}] \}_L ] =  \frac{1}{2} \{ 2^{n} - 2^{-n} \} = 0, \frac{3}{4}, \frac{15}{8}, \frac{63}{16}, \cdots.
  \label{sinh_seq}
  \end{align}
\end{definition}
Figure \ref{seq_coshx_sinhx}(a) shows difference between the usual hyperbolic cosine function ($\cosh x$) and the sequence version of the hyperbolic cosine function ($\{a_n[\cosh x]\}$) and Fig.~\ref{seq_coshx_sinhx}(b) shows difference between the usual hyperbolic sine function ($\sinh x$) and the sequence version of the hyperbolic sine function ($\{a_n[\sinh x]\}$).

The hyperbolic cosine function and hyperbolic sine function have the following properties.
  \begin{align}
    \frac{d}{dx} \cosh x &=  \sinh x,
  \label{cosh_prop}
  \\
    \frac{d}{dx} \sinh x &= \cosh x.    
  \label{sinh_prop}
  \end{align}
Then, analogously, we find that $\{ a_n [\cosh x] \}$ and $\{ a_n [ \sinh x] \}$ have the following properties.
  \begin{align}
    \mathcal{D}_R \{ a_n [\cosh x] \} &= \{ (2^{2n+1}-1)/2^{n+2} \} = \frac{1}{4}, \frac{7}{8}, \frac{31}{16}, \cdots =  \{ a_n [\sinh x] \} + \{ 1/(2^{n+2}) \},
  \\
    \mathcal{D}_R \{ a_n [\sinh x] \} &= \{ (2^{2n+1}+1)/2^{n+2} \} = \frac{3}{4}, \frac{9}{8}, \frac{33}{16}, \cdots =  \{ a_n [\cosh x] \} - \{ 1/(2^{n+2}) \}.
  \end{align}

Also, the following relation holds for the hyperbolic functions.
  \begin{align}
    \cosh^2x -  \sinh^2x = 1.  
  \label{cosh_sinh_prop}
  \end{align}
Analogously, the following relation holds for $\{ a_n [\cosh x] \}$ and $\{ a_n [ \sinh x] \}$.
  \begin{align}
    \{ a_n [ \cosh x] \}^2 - \{ a_n [ \sinh x] \}^2 = \{1\}. 
  \end{align}
  
Furthermore, there exists an addition theorem for the hyperbolic functions as follows.
  \begin{align}
    \cosh^2x + \sinh^2x =  \cosh 2x,
  \label{addition_hyper_1}
  \\
    2 \cosh x  \sinh x = \sinh 2x.
  \label{addition_hyper_2}
  \end{align}
Analogously, the following relation holds for $\{ a_n [\cosh x] \}$ and $\{ a_n [ \sinh x] \}$.
\begin{proposition} [Sequence version of addition theorem for hyperbolic function]
  \begin{align}
    \{ a_n [\cosh x] \}^2 + \{ a_n [\sinh x] \}^2 = \{ a_{2n} [ \cosh x] \},
  \\
    2 \{ a_n [\cosh x] \} \{ a_n [\sinh x] \} = \{ a_{2n} [\sinh x] \}.
  \end{align}
\end{proposition}

A sequence version of hyperbolic function can also be defined as follows by using Eq.~(\ref{ex-_natural}). 
\begin{definition}  [Sequence version of hyperbolic function (natural)]
  \begin{align}
    \{ a_n [\cosh x] \}^{\rm nat} = \frac{1}{2} [\{ a_n [[e^{x}] ] \} + \{ a_n [e^{-x}] \}^{\rm nat} ] = 1, 1, 2, 4, 8, \cdots,
  \\
    \{ a_n [\sinh x] \}^{\rm nat} = \frac{1}{2} [\{ a_n [[e^{x}] ] \} - \{ a_n [e^{-x}] \}^{\rm nat} ] = 0, 1, 2, 4, 8, \cdots.    
  \end{align}
\end{definition}
As an analogy of Eqs.~(\ref{cosh_prop}) and (\ref{sinh_prop}), $\{ a_n [\cosh x] \}^{\rm nat}$ and $\{ a_n [ \sinh x] \}^{\rm nat}$ have the following properties.
  \begin{align}
    \mathcal{D}_R \{ a_n [\cosh x] \}^{\rm nat} &= \{ a_n [\sinh x] \}^{\rm nat} ,
  \\
    \mathcal{D}_R \{ a_n [ \sinh x] \}^{\rm nat} &=  \{ a_n [\cosh x] \}^{\rm nat} .    
  \end{align}
Also, $\{ a_n [\cos x] \}^{\rm nat}$ and $\{ a_n [\sin x] \}^{\rm nat}$ have the following property analogous to Eq.~(\ref{cosh_sinh_prop}).
  \begin{align}
    (\{ a_n [\cosh x ] \}^{\rm nat})^2 - (\{ a_n [\sinh x ] \}^{\rm nat})^2= \{ a_n [e^{-x}]\}^{\rm nat}.
  \end{align}
Furthermore, the $ \{ a_n [\cosh x] \}^{\rm nat}$ and $ \{ a_n [ \sinh x] \}^{\rm nat}$ satisfy the following analogy of the addition theorem (Eq.~(\ref{addition_hyper_1}) and (\ref{addition_hyper_2})).
\begin{proposition} [Sequence version of addition theorem for hyperbolic function (natural)]
  \begin{align}
    (\{ a_n [ \cosh x] \}^{\rm nat})^2 + (\{ a_n [ \sinh x] \}^{\rm nat})^2 = \{ a_{2n} [\cosh x] \}^{\rm nat},
  \\
    2 \{ a_n [ \cosh x] \}^{\rm nat} \{ a_n [ \sinh x] \}^{\rm nat} = \{ a_{2n} [ \sinh x] \}^{\rm nat}.
  \end{align}
\end{proposition}

Figure \ref{seq_coshx_sinhx}(c) shows difference between $\{a_n[\cosh x]\}$ ($\{a_n[\sinh x]\}^{\rm nat}$) and $\{a_n[\cosh x]\}^{\rm nat}$ ($\{a_n[\sinh x]\}$).

\newpage
\section{Trigonometric function}\label{Sec:trigonometric}
\begin{figure}[h]
 \centering
 \includegraphics[width=14.5cm]{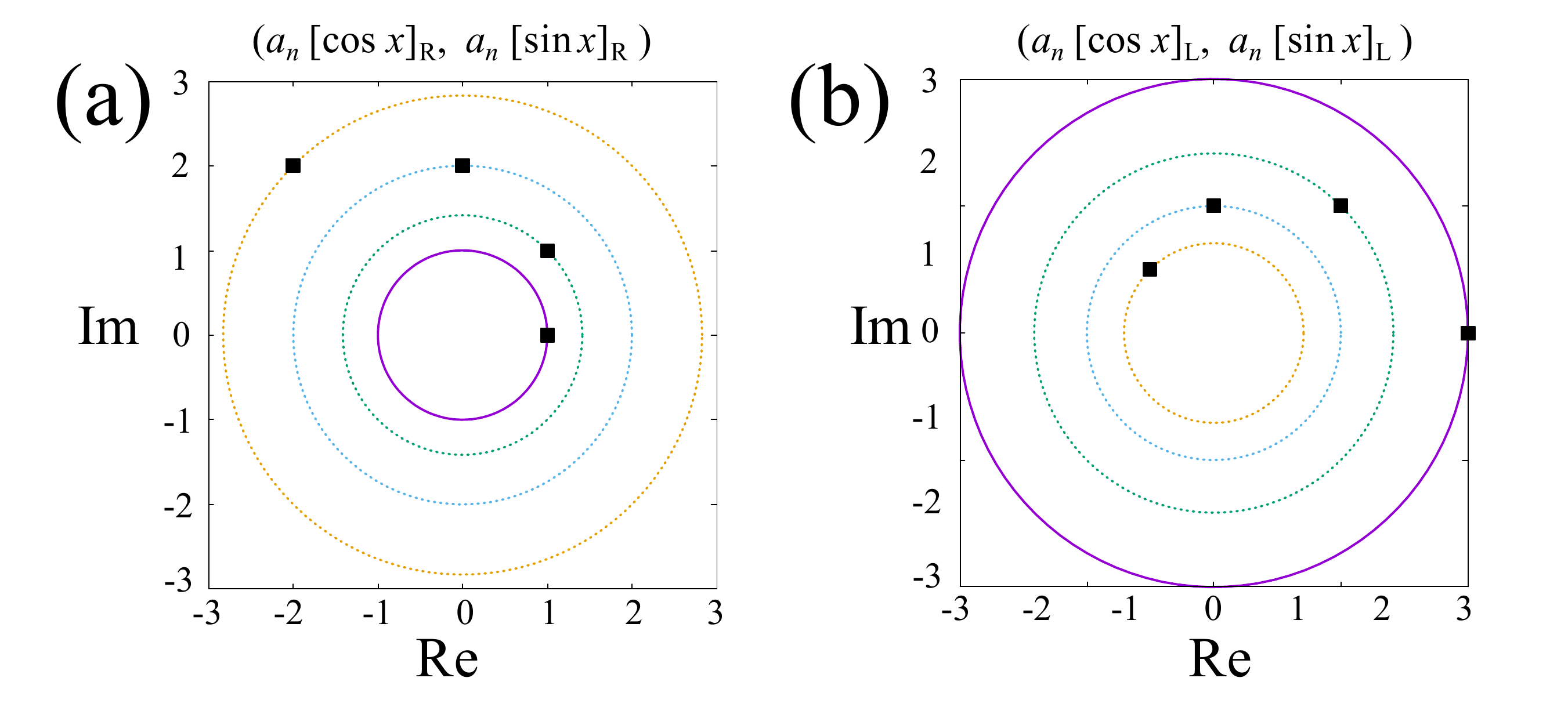}
 \caption{Graph of (a) $( \{ a_n [ \cos x] \}_R,  \{ a_n [ \sin x] \}_R)$ corresponding to the divergent spiral and (b) $( \{ a_n [ \cos x] \}_L,  \{ a_n [ \sin x] \}_L)$ corresponding to the convergent spiral.}
 \label{seq_cosx_sinx}
\end{figure}

First, we define sequence version of trigonometric function in two ways. By extending domain of Eq.~(\ref{def_ex_gene1}) from real numbers to complex numbers and focusing on the following function relation
  \begin{align}
    e^{ix} = \cos x + i \sin x,
  \end{align}
we obtain the following sequence version of the cosine function and sin function.
\begin{definition} [Right sequence version of trigonometric function]
  \begin{align}
    \{ a_n [ \cos x] \}_R &= \{ {\rm Re} (1+i)^n \} = 1,1,0,-2,-4,-4,0,8,16,16,\cdots,
  \label{Yuki1}
  \\
    \{ a_n [ \sin x] \}_R &= \{ {\rm Im} (1+i)^n \} = 0,1,2,2,0,-4,-8,-8,0,16,\cdots.
  \label{Yuki2}
  \end{align}
\end{definition}
Figure \ref{seq_cosx_sinx}(a) shows a graph of $( \{ a_n [ \cos x] \}_R,  \{ a_n [ \sin x] \}_R)$ that represents the divergent spiral. This contrasts with the circle represented by $(\cos x, \sin x)$.

The cosine function and sine function have the following properties for differential.
  \begin{align}
    \frac{d}{dx}  \cos x &= -  \sin x,
  \label{cs_1_1}
  \\
    \frac{d}{dx}  \sin x &=  \cos x.   
  \label{cs_1_2}
  \end{align}
Then, analogously, we find that $\{ a_n [{\rm cos}\,x] \}_R $ and $\{ a_n [ \sin x] \}_R$ have the following properties.
  \begin{align}
    \mathcal{D}_R \{ a_n [\cos x] \}_R &= 0,-1,-2,-2,0,4,8,8,0,-16,\cdots = - \{ a_n [\sin x] \}_R ,
  \\
    \mathcal{D}_R \{ a_n [ \sin x] \}_R &= 1,1,0,-2,-4,-4,0,8,16,16, \cdots = \{ a_n [ \cos x ] \}_R.  
  \end{align}

Also, the following relation holds for the trigonometric functions.
  \begin{align}
      \cos^2x + \sin^2x = 1.
  \label{cs_2}
  \end{align}
Analogously, the following formula for $\{ a_n [ \cos x] \}_R $ and $\{ a_n [ \sin x] \}_R$ holds. 
  \begin{align}
     \{ a_n [ \cos x] \}_R^2 +  \{ a_n [ \sin x] \}_R^2 = 2^n.
  \end{align}
In addition, the tangent function is given by
  \begin{align}
    \tan x = \frac{ \sin x}{ \cos x}.
  \end{align}
Then, the sequence version of the tangent function determined by $\{ a_n [ \cos x] \}_R $ and $\{ a_n [ \sin x] \}_R$ is given by
  \begin{align}
    \{ a_n [\tan x] \}_R = \frac{\{ a_n [ \sin x] \}_R}{\{ a_n [ \cos x] \}_R } = 0,1,\infty,-1,0, 1,-\infty,-1,0,1,\cdots.
  \end{align}
where we wrote $n/0=\infty$ and $-n/0=-\infty$ for a natural number $n$.

Furthermore, the trigonometric functions satisfy an addition theorem as follows.
  \begin{align}
    \cos 2x &= \cos^2x - \sin^2x,
  \label{cs_3}
  \\
    \sin 2x &= 2 \sin x  \cos x.
  \label{cs_4} 
  \end{align}
Then, we find the direct analogy of the addition theorem for trigonometric function as follows.
\begin{proposition} [Addition theorem for right sequence version of trigonometric function]
  \begin{align}
    \{ a_n [\cos x] \}_R^2 - \{ a_n [\sin x] \}_R^2 
      &= 1,0,-4,0,16,0,\cdots = \{ a_{2n} [\cos x] \}_R,
  \\
    2  \{ a_n [\cos x] \}_R  \{ a_n [\sin x] \}_R
      &= 0,2,0,-8,0,32,\cdots = \{ a_{2n} [\sin x] \}_R.
  \end{align}
\end{proposition}

By using Eq.~(\ref{def_ex_gene2}), one has another definition of sequence version of cosine function and sine function as follows.
\begin{definition} [Left sequence version of trigonometric function]
  \begin{align}
    \{ a_n [ \cos x] \}_L &= \{ {\rm Re} (1-i)^{-n} \} = 1,\frac{1}{2}, 0, -\frac{1}{4}, -\frac{1}{4}, -\frac{1}{8}, 0, \frac{1}{16} \frac{1}{16}, \frac{1}{32}, \cdots,
  \\
    \{ a_n [\sin x] \}_L &= \{ {\rm Im} (1-i)^{-n} \} = 0,\frac{1}{2}, \frac{1}{2}, \frac{1}{4}, 0, -\frac{1}{8}, -\frac{1}{8}, -\frac{1}{16}, 0, \frac{1}{32}, \cdots.
    \end{align}
\end{definition}
Figure \ref{seq_cosx_sinx}(b) shows a graph of $( \{ a_n [ \cos x] \}_L,  \{ a_n [ \sin x] \}_L)$ that represents the convergent spiral. This contrasts with the circle represented by $(\cos x, \sin x)$, and also the divergent spiral represented by $( \{ a_n [ \cos x] \}_R,  \{ a_n [ \sin x] \}_R)$.

As an analogy of Eqs.~(\ref{cs_1_1}) and (\ref{cs_1_2}), $\{ a_n [\cos x] \}_L $ and $\{ a_n [\sin x] \}_L$ have the following properties.
  \begin{align}
    \mathcal{D}_L \{ a_n [\cos x] \}_L &=  0, -\frac{1}{2}, -\frac{1}{2}, -\frac{1}{4}, 0, \frac{1}{8}, \frac{1}{8}, \frac{1}{16}, 0, -\frac{1}{32}, \cdots = - \{ a_n [\sin x] \}_L ,
  \\
    \mathcal{D}_L \{ a_n [\sin x] \}_L &= 1,\frac{1}{2}, 0, -\frac{1}{4}, -\frac{1}{4}, -\frac{1}{8}, 0, \frac{1}{16} \frac{1}{16}, \frac{1}{32}, \cdots = \{ a_n [\cos x] \}_L, 
  \end{align}
where we assumed $a_{-1} [\cos x]_L = 1$ and $a_{-1} [\sin x]_L = -1$.

Also, $\{ a_n [\cos x] \}_L$ and $\{ a_n [\sin x] \}_L$ have the following property that is analogous to Eq.~(\ref{cs_2}). 
  \begin{align}
     \{ a_n [\cos x] \}_L^2 +  \{ a_n [\sin x] \}_L^2 = 2^{-n}.
  \end{align}

In addition, the left sequence version of the tangent function is given by
  \begin{align}
    \{ a_n [\tan x] \}_L = \frac{\{ a_n [\sin x] \}_L}{\{ a_n [\cos x] \}_L } = 0,1,\infty,-1,0, 1,-\infty,-1,0,1,\cdots = \{ a_n [\tan x] \}_R.
  \end{align}

Furthermore, $ \{ a_n [\cos x] \}_L$ and $ \{ a_n [ \sinh x] \}_L$ satisfy the following formula that is analogous to the addition theorem (Eq.~(\ref{cs_3}) and (\ref{cs_4})).
\begin{proposition} [Addition theorem for left sequence version of trigonometric function]
  \begin{align} 
    \{ a_n [ \cos x] \}_L^2 - \{ a_n [ \sin x] \}_L^2 
      &= 1,0,-\frac{1}{4},0,\frac{1}{16},0,\cdots = \{ a_{2n} [ \cos x] \}_L,
  \\
    2  \{ a_n [\cos x] \}_L  \{ a_n [\sin x] \}_L
      &= 0,\frac{1}{2},0,-\frac{1}{8},0,\frac{1}{32},\cdots = \{ a_{2n} [\sin x] \}_L.
  \end{align}
\end{proposition}

\begin{figure}[h]
 \centering
 \includegraphics[width=14.5cm]{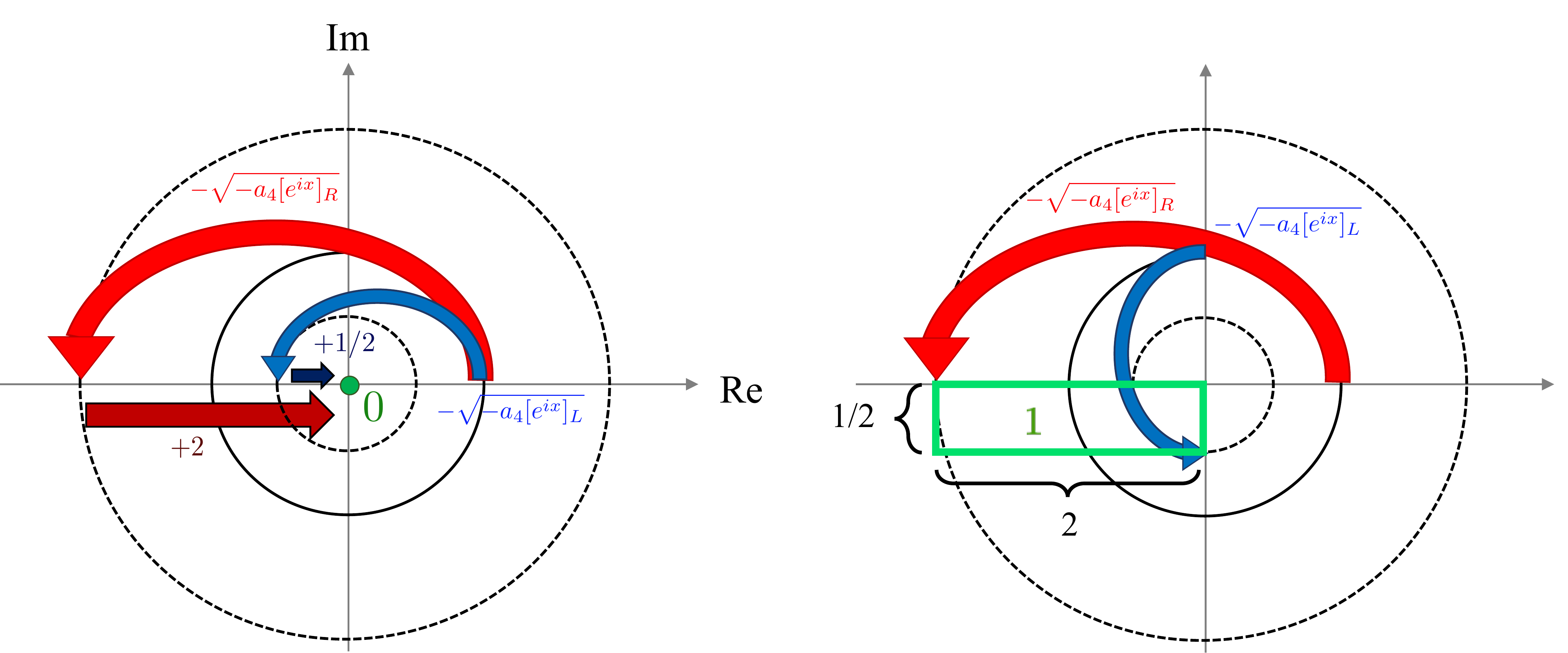}
 \caption{(a) Sequence version of Euler's identity ($a_4 [e^{ix}]_R + 2^2 = 0$ and $a_4 [e^{ix}]_L + 2^{-2} = 0$), (b) $a_4 [e^{ix}]_L a_4 [e^{ix}]_R = 1.$}
 \label{seq_Euler}
\end{figure}

The following equality is called Euler's identity and connects the Napier's constant ($e$), the additive identity ($0$), the multiplicative identity ($1$), the imaginary unit ($i$), the circle ratio ($\pi$) , and would be one of the most famous and beautiful equality in mathematics. 
  \begin{align}
    e^{i\pi} + 1 = 0.
  \end{align}
Note that this equation can be rewritten in two ways as follows.
  \begin{align}
      e^{i\pi} + 1^2 = 0,
  \label{Euler_orig1}
  \\
      e^{i\pi} + 1^{-2} = 0.
  \label{Euler_orig2} 
  \end{align}
Then, we construct sequence version of the Euler's identity by using Eq.~(\ref{def_ex_gene1}) and Eq.~(\ref{def_ex_gene2}) as follows.
\begin{theorem} [Sequence version of Euler's identity] 
  \begin{align}
     a_4 [e^{ix}]_R + 2^2 = 0,
  \label{Euler_seq1}
  \\
     a_4 [e^{ix}]_L + 2^{-2} = 0. 
  \label{Euler_seq2}
  \end{align}
\end{theorem}
Figure \ref{seq_Euler}(a) shows schematic representation of Eqs.~(\ref{Euler_seq1}) and (\ref{Euler_seq2}).

Equation (\ref{Euler_seq1}) is obtained by replacing $e^{i\pi}$ by $a_4 [e^{ix}]_R$ and $1^2$ by $2^2$ in Eq.~(\ref{Euler_orig1}). Also, Eq.~(\ref{Euler_seq2}) is obtained by replacing $e^{i\pi}$ by $a_4 [e^{ix}]_L$ and $1^{-2}$ by $2^{-2}$ in Eq.~(\ref{Euler_orig2}). In addition, as an analogy of the relation (see Fig.~\ref{Euler}(c))
  \begin{align}
    e^{i\pi} e^{i\pi} = 1,
  \end{align}
the following formula holds.
  \begin{align}
    a_4 [e^{ix}]_L a_4 [e^{ix}]_R = 1.
  \label{Euler_seq3}
  \end{align}
Figure \ref{seq_Euler}(b) shows schematic representation of Eq.~(\ref{Euler_seq3}).

\begin{figure}[h]
 \centering
 \includegraphics[width=14.5cm]{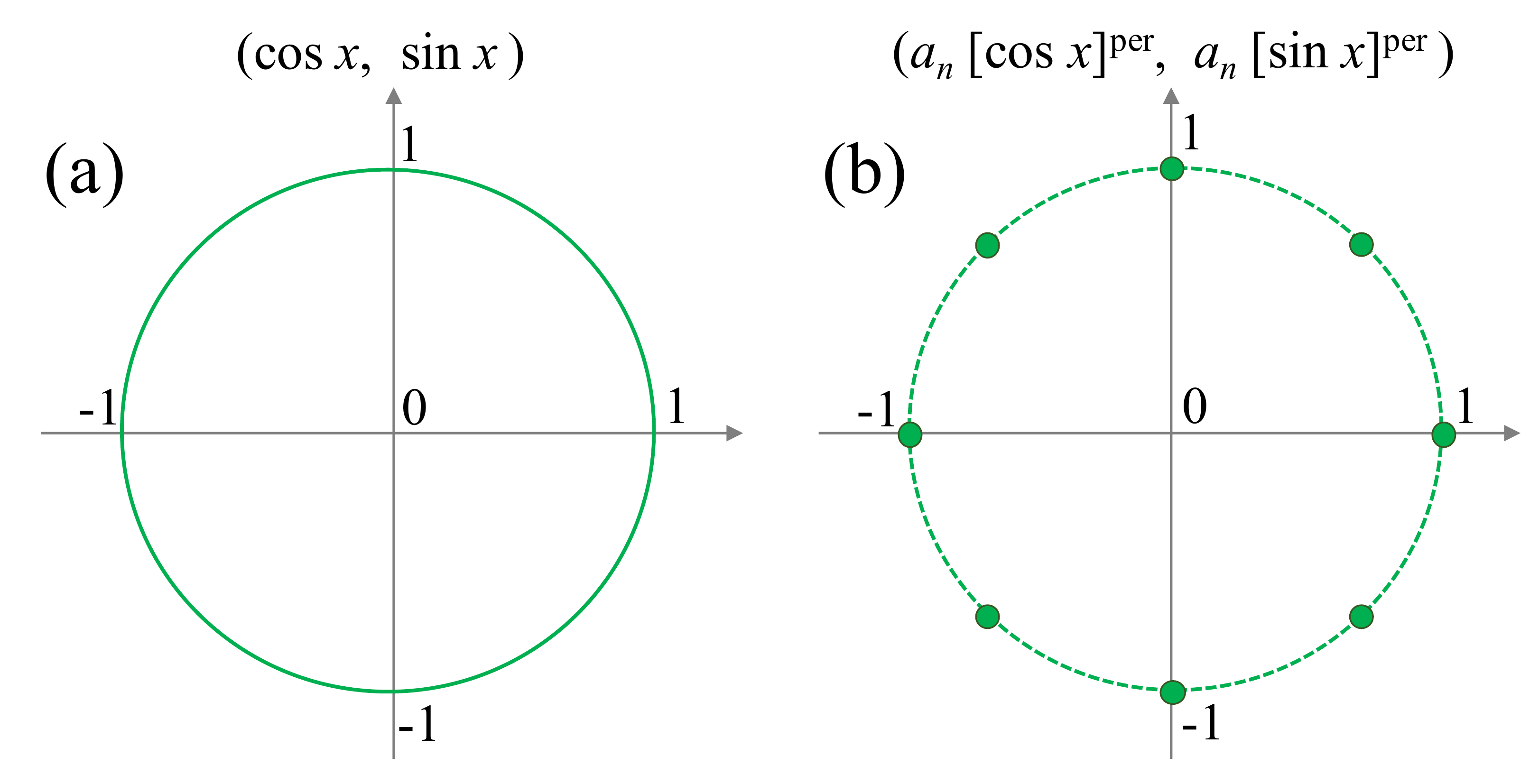}
 \caption{(a) Representation of $(\cos x, \sin x)$ by the unit circle. (b) Representation of $(a_n[\cos x]^{\rm per}, a_n[\sin x]^{\rm per})$ by the eight points on the unit circle.}
 \label{seq_cosx_sinx_per}
\end{figure}

Also, by combining $\{ a_n [ \cos x] \}_R$ ($\{ a_n [ \sin x] \}_L$) and $\{ a_n [ \cos x] \}_L$ ($\{ a_n [ \sin x] \}_L$) generally, one has the following sequence version of $\cos x$ and $\sin x$ whose periodicity is 8. 
  \begin{align}
    \{ a_n [\cos x] \}^{\rm per} &= \{ {\rm sgn} \ (a_n[\cos x]_R ) \} \times  \sqrt{ \{ a_n [ \cos x] \}_R \times \{ a_n [ \cos x] \}_L} 
  \\
    &= 1, \frac{1}{\sqrt{2}}, 0, -\frac{1}{\sqrt{2}}, -1, -\frac{1}{\sqrt{2}}, 0, \frac{1}{\sqrt{2}}, \cdots,
  \\
    \{ a_n [\sin x] \}^{\rm per} &= \{ {\rm sgn} \ (a_n[\sin x]_R ) \} \times \sqrt{ \{ a_n [ \sin x] \}_R \times \{ a_n [ \sin x] \}_L} 
  \\
    &= 0, \frac{1}{\sqrt{2}}, 1, \frac{1}{\sqrt{2}}, 0, -\frac{1}{\sqrt{2}}, -1, -\frac{1}{\sqrt{2}}, \cdots.
  \end{align}
This periodicity is analogous to $2\pi$-periodicity of $\cos x$ and $\sin x$. In Fig.~\ref{seq_cosx_sinx_per}, we show representation of $(\cos x, \sin x)$ by the unit circle and of $(a_n[\cos x]^{\rm per}, a_n[\sin x]^{\rm per})$ by the eight points on the unit circle. While the $2\pi$ is the periodicity of $\cos x$ and $\sin x$ and the circumference of a unit circle, and the 8 is the periodicity of $\{a_n[\cos x]\}^{\rm per}$ and $\{a_n[\sin x]\}^{\rm per}$ and the length of the circumference of the unit circle when the unit is given by the distance between two adjacent points on the unit circle ($\pi/4$).

\newpage
\section{Fibonacci sequences}\label{Sec:fibonacci}
In this section, we apply the sequence version of Maclaurin series to the Fibonacci sequence and its generalizations.

\subsection{Fibonacci sequence}
\begin{figure}[h]
 \centering
 \includegraphics[width=14.5cm]{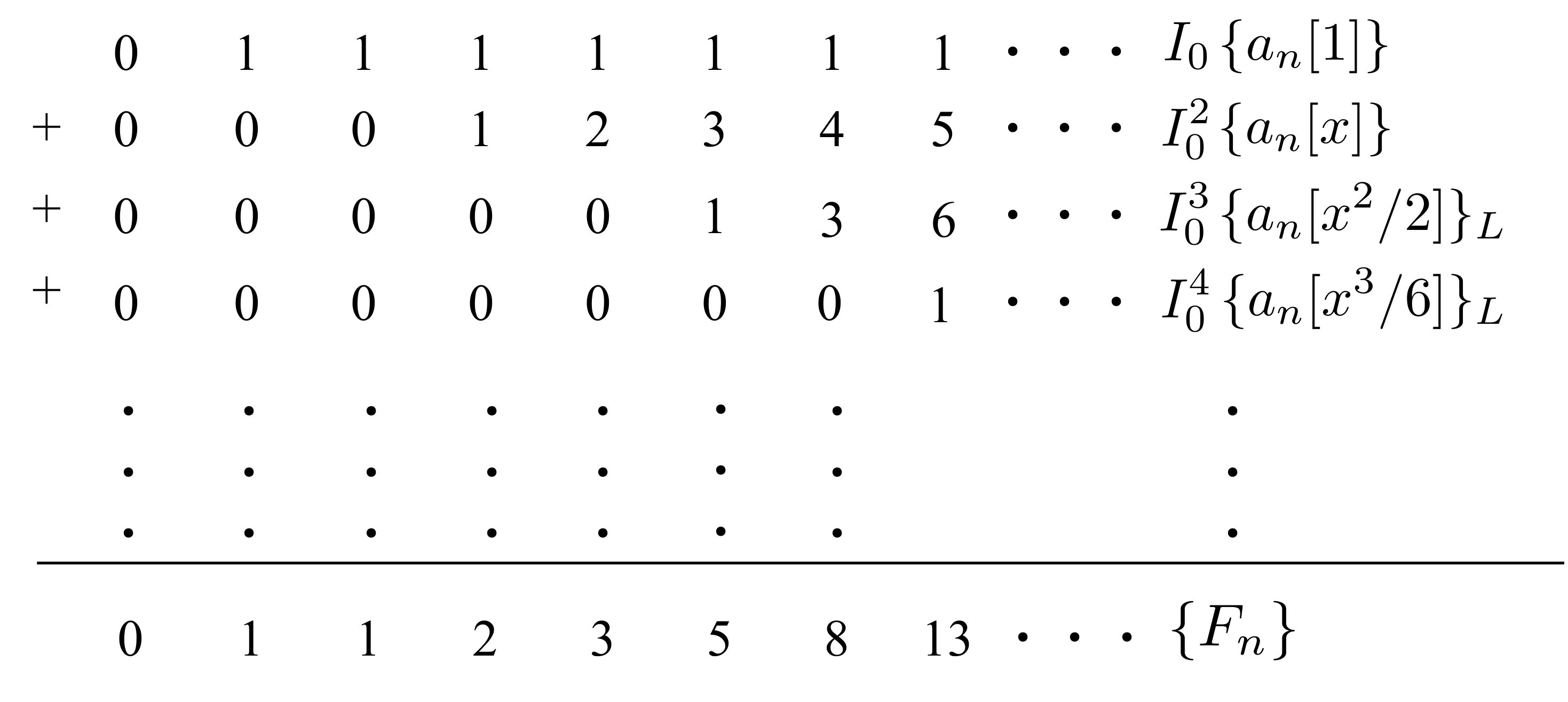}
 \caption{Digrammatic representation of sequence version of Maclaurin series for Fibonacci sequence corresponding to Pascal's triangle.}
 \label{Fibo}
\end{figure}

\begin{figure}[h]
 \centering
 \includegraphics[width=14.5cm]{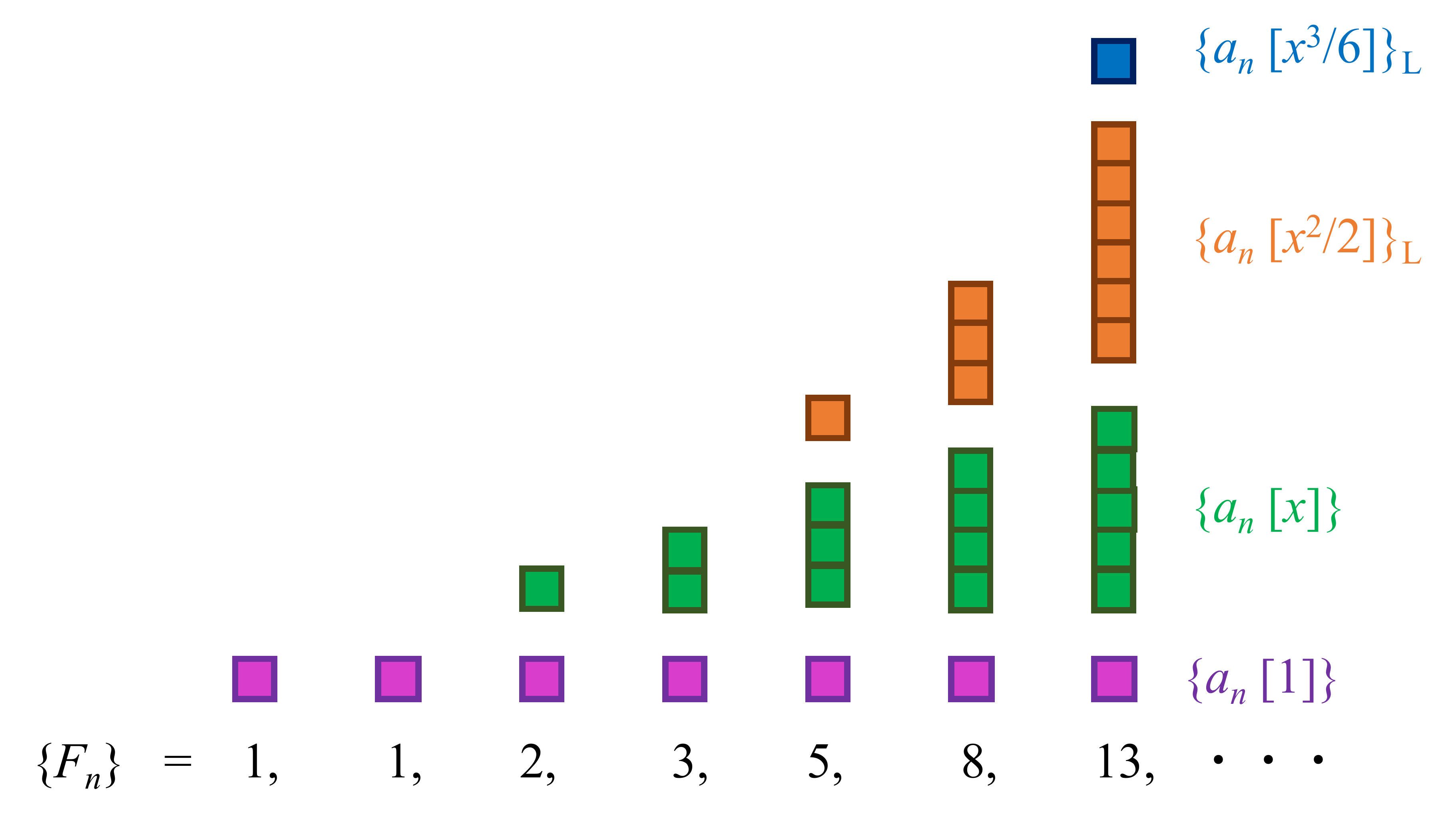}
 \caption{Diagrammatic representation of sequence version of Maclaurin series for Fibonacci sequence by using blocks.}
 \label{block_Fibonacci}
\end{figure}

Fibonacci sequence is defined by
  \begin{align}
    \{ F_n \} = 0, 1, 1, 2, 3, 5, 8, \cdots,
  \end{align}
and satisfies the relation
  \begin{align}
   \lim_{n\to\infty} \frac{ F_{n}  }{ F_{n-1}  } = \phi,
  \end{align}
where $\phi$ is the golden ratio and also the relation
  \begin{align}
     F_{n+2}   =  F_{n+1}  +  F_n .
  \end{align}
This relation can be rewritten by using the right differential $\mathcal{D}_R$ defined in Eq.~(\ref{DR}) as follows.
  \begin{align}
    \mathcal{D}_R \{ F_{n+1} \} = \{ F_n \} .
  \end{align}
Also, by using the shift operator $\mathcal{S}_k$ defined in Eq.~(\ref{shift}), one has 
  \begin{align}
    \mathcal{D}_R \mathcal{S}_1 \{ F_n \} = \{ F_n \} .  
  \end{align}
This is similar to Eq.~(\ref{def_ex_DR}) that the right sequence version of exponential function satisfies. Then, we obtain the following theorem by using the insertion operator $I_a$ defined in Eq.~(\ref{insertion}) that is similar to Eq.~(\ref{exponential2}) (see Fig.~\ref{Fibo} and compare Fig.~\ref{block_Fibonacci} with Fig.~\ref{block_ex}).
\begin{theorem} [Sequence version of Maclaurin series for Fibonacci sequence]
  \begin{align}
    \{ F_n \} = I_0 \{ a_n [1] \} + I_0^2 \{ a_n [x] \} + I_0^3 \{ a_n [x^2/2] \} + \cdots = \sum_{n=0}^{\infty} I_0^{n+1} \{ a_n [x^n/n!] \}.
  \end{align}
\end{theorem}

Also, the negafibonacci sequence \cite{negaFibonacci} is defined by
  \begin{align}
    \{ F_n^- \} = 0, 1, -1, 2, -3, 5, -8, \cdots = \{ F_n \} \times \{ (-1)^{n+1} \},
  \end{align}
and satisfies the following relation.
  \begin{align}
    \mathcal{D}_R  \{ F_n^- \} = - \mathcal{S}_2\{  F_n^- \}.
  \end{align}
Although this is similar to $e^{-x}$, however the relation $\{ F_n \}  \times \{ F_n^- \} = \{ a_n [1]\}$ does not hold. Instead the following relation holds that is equivalent to Cassini's Fibonacci identity \cite{Cassini}.
  \begin{align}
    \{ F_n \} \times \mathcal{S}_2 \{ F_n^- \} + \mathcal{S}_1 \{ F_n \} \times \mathcal{S}_1 \{ F_n^- \} = \{ a_n [1] \}.
  \end{align}
For example, $F_1 F_3^- + F_2^-F_2 = 1\times2 + 1\times(-1) = 1$.

\subsection{$(P,\,Q)$-Fibonacci sequence}
\begin{figure}[h]
 \centering
 \includegraphics[width=14.5cm]{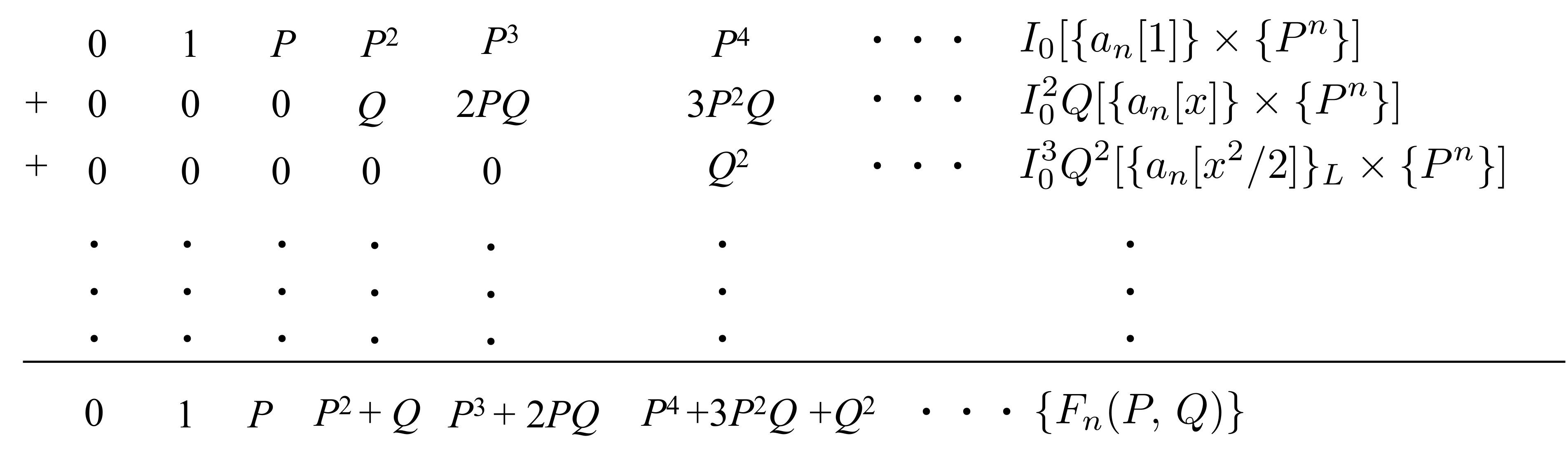}
 \caption{Digrammatic representation of sequence version of Maclaurin series for $(P, \ Q)$-Fibonacci sequence.}
 \label{PQ_Fibo}
\end{figure}

The $(P,\,Q)$-Fibonacci sequence is given by
  \begin{align}
    \{ F_n (P, Q)\} = 0, 1, P, P^2+Q, P^3+2PQ, P^4+3P^2Q+Q^2,\cdots,
  \end{align}
and is a generalization of Fibonacci sequence and satisfies the following relation \cite{gene:Fibo, gene:Fibo2}.
  \begin{align}
    \{ F_{n+2} (P, Q)\} = P \{ F_{n+1} (P, Q)\} + Q \{ F_{n} (P, Q)\} .
  \end{align}
We find the $(P,\,Q)$-Fibonacci sequence has the following sequence version of Maclaurin series (Fig.~\ref{PQ_Fibo}).
\begin{proposition} [Sequence version of Maclaurin series for $(P,Q)$-Fibonacci sequence]
  \begin{align}
    \{ F_{n} (P, Q)\} &= I_0 [\{ a_n [1] \} \times  \{ P^n \} ] + Q I_0^2 [ \{ a_n [x] \} \times  \{ P^n \} ]  +  \cdots,
  \\
    &= \sum_{k=0}^{\infty} I_0^{k+1} Q^k [\{ a_n [x^k/k!] \} \times  \{ P^n \} ]    .
  \end{align}
\end{proposition}

\subsubsection{Example 1: Pell sequence ($P=2$, $Q=1$)}
\begin{figure}[h]
 \centering
 \includegraphics[width=14.5cm]{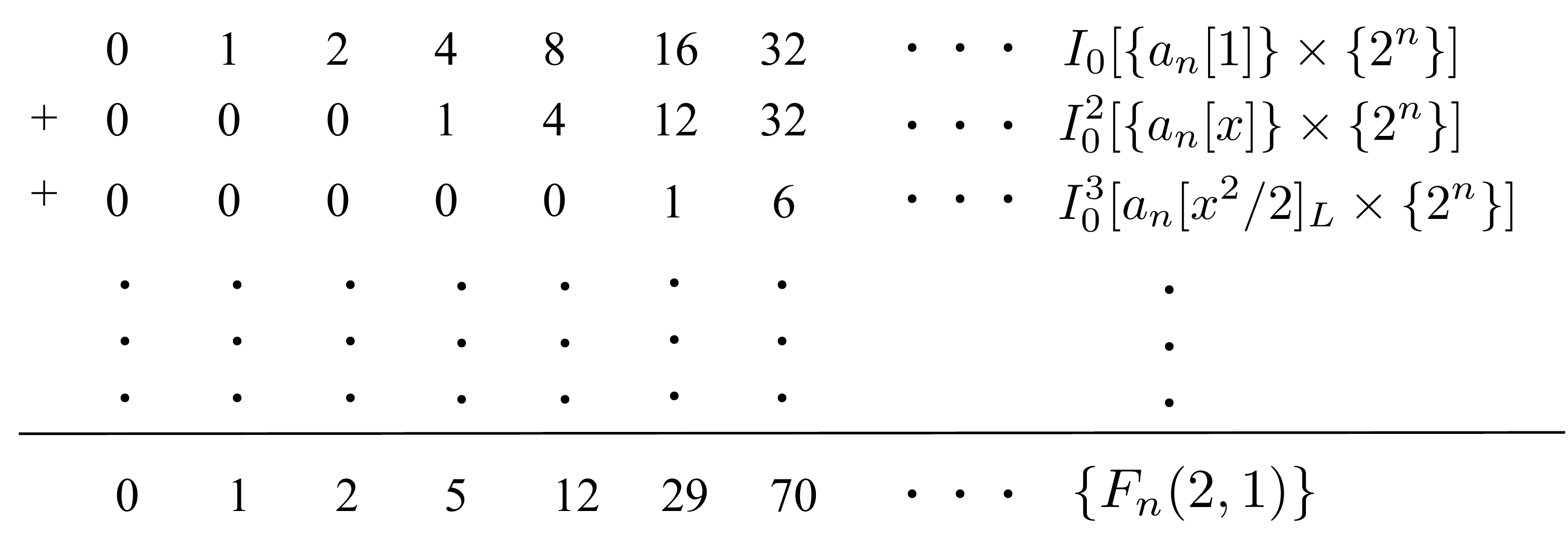}
 \caption{Digrammatic representation of sequence version of Maclaurin series for $(2, \ 1)$-Fibonacci sequence (Pell sequence).}
 \label{Pell}
\end{figure}

The $(2,1)$-Fibonacci sequence is called Pell sequence \cite{Pell, Pell2, Pell3, Pell4} and satisfies the following relation.
  \begin{align}
   \lim_{n\to\infty} \frac{ F_{n} (2, 1) }{ F_{n-1} (2, 1) } = 1 + \sqrt{2},
  \end{align}
where $1 + \sqrt{2}$ is the silver ratio. The sequence version of Maclaurin series for Pell Sequence is shown in Fig.~\ref{Pell}.

\subsubsection{Example 2: Jacobsthal sequence ($P=1$, $Q=2$)}
\begin{figure}[h]
 \centering
 \includegraphics[width=14.5cm]{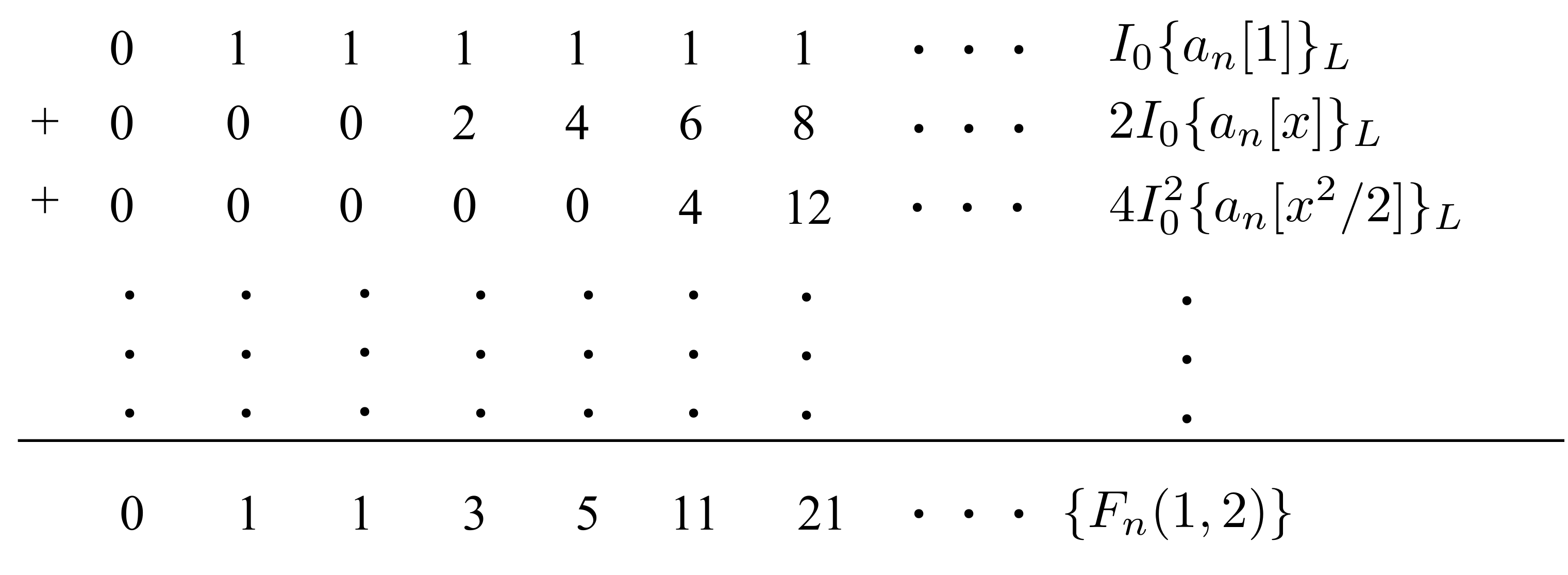}
 \caption{Digrammatic representation of sequence version of Maclaurin series for $(1, \ 2)$-Fibonacci sequence (Jacobsthal sequence).}
 \label{Jacobsthal}
\end{figure}

The $(1,2)$-Fibonacci sequence is called Jacobsthal sequence \cite{Jacobsthal, Jacobsthal2, Jacobsthal3} and the sequence version of Maclaurin series for Jacobsthal Sequence is shown in Fig.~\ref{Jacobsthal}.

\subsection{$k$-bonacci sequence}
In this subsection, we propose the sequence version of $k$-bonacci sequence \cite{k-bonacci} satisfying the following relation.
  \begin{align}
    \{ F_{n+k} (k) \} = \sum_{l=0}^{k-1} \{ F_{n+l} (k) \}.
  \end{align}
Then, we construct the sequence version of Maclaurin series for $ \{ F_{n+k} (k) \}$. First, we define deformed integral by combining the integral $\mathcal{I}_0$ and the insertion $I_0$.
\begin{definition} [Deformed integral]
For an integer $k\geq2$,
  \begin{align}
    \tilde{\mathcal{I}}_k = \sum_{l=0}^{k-2} I_0^l \mathcal{I}_L^0.
  \end{align}
\end{definition}
Then, we have the following proposition.
\begin{proposition} [Sequence version of Maclaurin series for $k$-bonacci sequence]
For an integer $k$ such that $k\geq2$,
  \begin{align}
     \{ F_{n} (k) \} &=  I_0^{k-1} \{ a_n [1] \} + \tilde{\mathcal{I}}_k I_0^{k} \{ a_n [1] \} + \tilde{\mathcal{I}}_k^2 I_0^{k+1} \{ a_n [1] \} + \cdots
  \\
    &= \sum_{l=0}^{\infty} (\tilde{\mathcal{I}}_k)^l I_0^{k-1+l} \{ a_n [1] \}.
  \end{align}
\end{proposition}
This is a generalization of the left sequence version of Maclaurin series for right sequence version of exponential function (Prop.~\ref{exponential}). Also, we define the limit of sequence as follows. 
\begin{definition} [Limit of sequence] For an integer $k$ and a sequence $\{ a_n \}$ depending on $k$,
  \begin{align}
    \lim_{k\to\infty} \{ a_n (k) \} =  \lim_{k\to\infty} a_0(k),  \lim_{k\to\infty} a_1(k), \cdots.
  \end{align}
\end{definition}
Then, we obtain the following relation that states $\mathcal{S}_{k} \{ F_{n} (k) \}$ approaches to $\{ a_n [e^x] \}_R$ in the limit of $k\to\infty$.
\begin{proposition} [Limit of $k$-bonacci sequence and right sequence version of exponential function]
  \begin{align}
    \lim_{k\to\infty} \mathcal{S}_{k} \{ F_{n} (k) \}  = \{ a_n [e^x] \}_R.
  \end{align}
\end{proposition}

\begin{figure}[h]
 \centering
 \includegraphics[width=14.5cm]{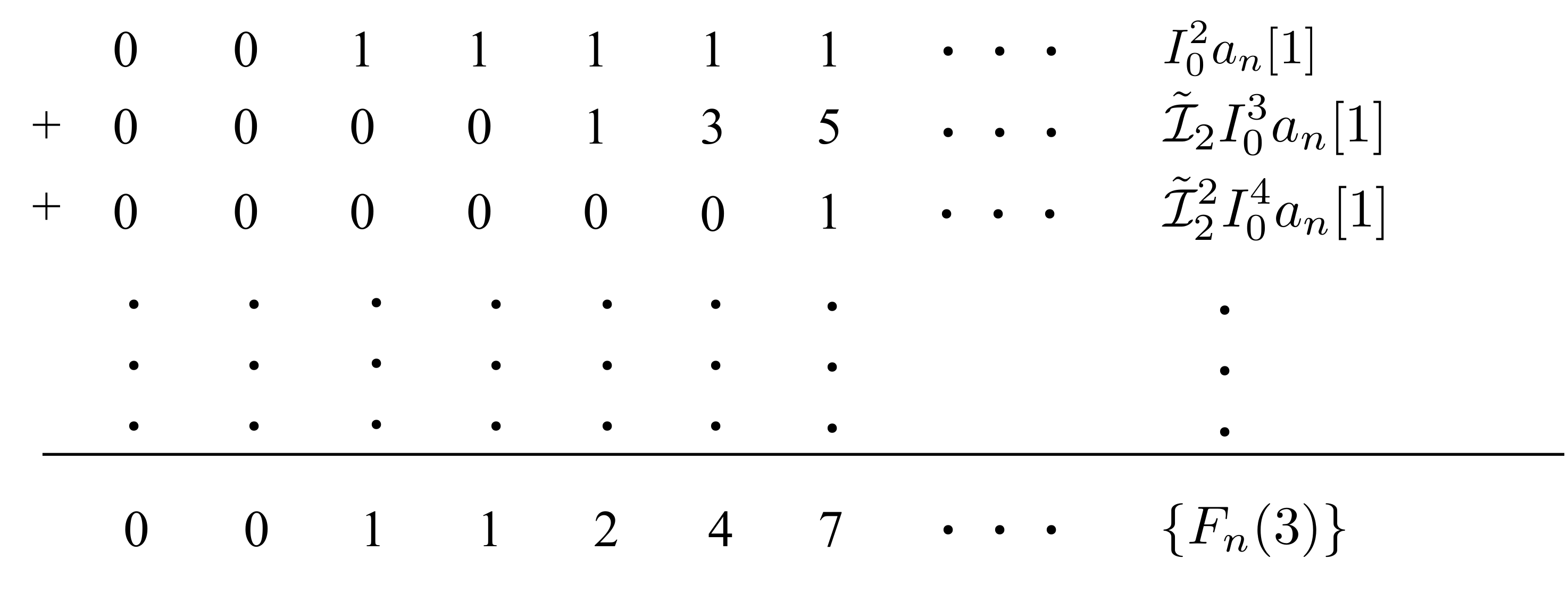}
 \caption{Digrammatic representation of sequence version of Maclaurin series for tribonacci sequence ($\{ F_{n} (3) \}$).}
 \label{tribonacci}
\end{figure}

We show sequence version of Maclaurin series for tribonacci sequence ($\{ F_{n} (3) \}$) in Fig.~\ref{tribonacci}. For example, the following relation holds.
  \begin{align}
    \tilde{\mathcal{I}}_2 I_0^3 \{ a_n [1] \} &= \mathcal{I}_L^0 [ I_0^3 \{ a_n [1] \} +  I_0^4 \{ a_n [1] \} ] = \mathcal{I}_L^0 (0, 0, 0,  1, 2, 2, 2, \cdots)
  \nonumber\\
    &= 0,0,0,0,1,3,5,7, \cdots.
  \end{align}

\newpage
\section{Factorial and Bell number}\label{Sec:factorial:Bell}
In this section, we construct sequence dual of factorials and Bell numbers and propose sequence dual of the modular property of factorial concerning prime number and of Bell numbers concerning prime number.

\subsection{Factorial}
Factorial is defined for a non-negative number $n$ as follows.
  \begin{align}
    n! = n \times (n-1) \times \cdots \times 2 \times 1,
  \end{align}
and $0!=1$. Because the factorial has the following integral representation
  \begin{align}
    n! = \int_0^{\infty} x^{n-1} e^{-x} dx,
  \end{align}
we construct sequence dual of factorial as follows.
  \begin{align}
    \{ a_n \langle n! \rangle \} &= \sum_{k=0}^{\infty} k^n \times  2^{-k} ,
  \\
    & = 2, 2, 6, 26, 150, 1082, 9366, 94586, 1091670, 14174522, \cdots.
  \end{align}
We find that $ \{ a_n \langle n! \rangle \} $ deeply relates to the factorial $n!$. First, we focus on the following transform for a sequence $\{ a_n \}$ called Stirling transform \cite{Stirling:transform}.
  \begin{align}
    \{ b_n \} =  \left\{ \sum_{k=0}^n S(n,k) a_k \right\},
  \end{align}
where $S(n,k)$ is the Stirling number of the second kind that counts how many ways to partition a set of $n$ elements into $k$ non-empty subsets. For example, $S(3,2)=3$ because $\{a, b, c\}$ has the following subgroups that has two elements: $\{a,b\}$, $\{b,c\}$, $\{c,a\}$.

Then, the following proposition holds.
\begin{proposition} [Sequence dual of factorial and Stirling transform]
  \begin{align}
    \mathcal{S}_1 \{ a_n \langle n! \rangle \}
       = 2 \left\{ \sum_{k=0}^n S(n,k) k! \right\}.
  \label{factorial_Stirling}
  \end{align}
\end{proposition}
This proposition shows the sequence dual of factorial is equivalent to the Stirling transform of the second kind for the factorial up to the shift operation and scalar multiplication.

\begin{figure}[h]
 \centering
 \includegraphics[width=10.0cm]{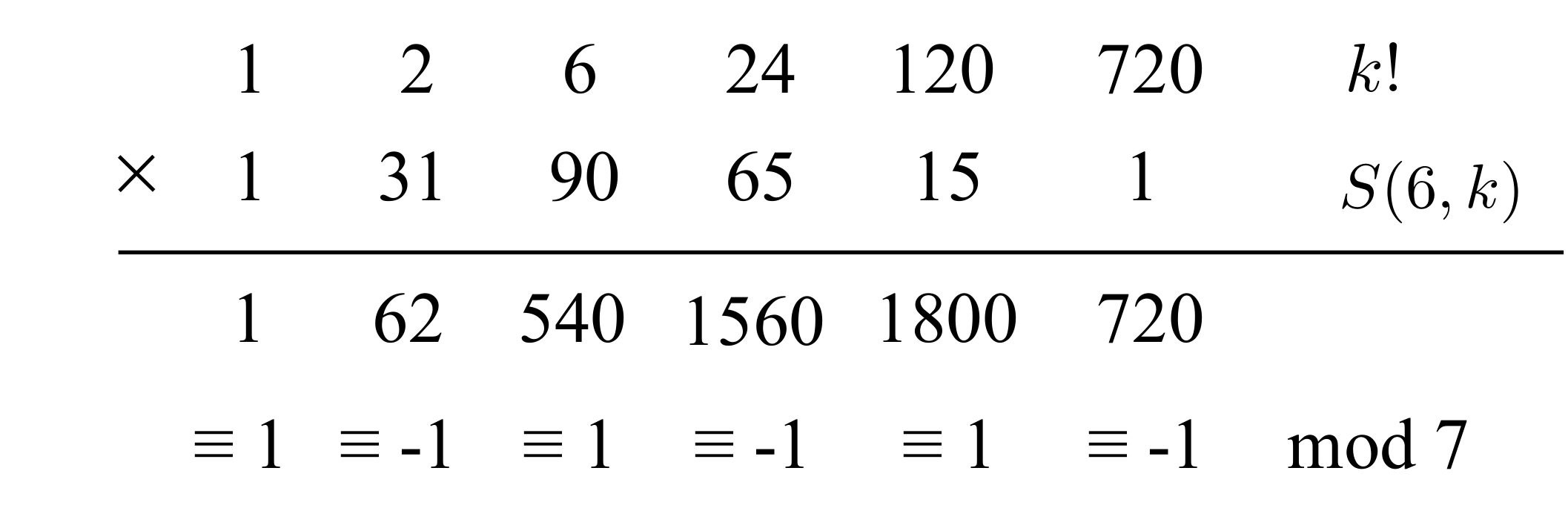}
 \caption{Diagrammatic representation of Eq.~(\ref{prime:factorial:2}) for $n=7$.}
 \label{Factorial_prime}
\end{figure}

Next, we focus on the following property of the factorial concerning the prime numbers called Wilson's theorem. For a natural number $n$,
  \begin{align}
    n! \equiv -1 \ ({\rm mod} \ n+1),
  \end{align}
if and only if $n+1$ is a prime. For example, $4\times3\times2\times1=24=5\times5-1$, $5\times4\times3\times2\times1=120=6\times20-0$. Then, we obtain the sequence dual of Wilson's theorem as follows.
\begin{theorem} [Sequence dual of Wilson's theorem]
For a natural number $n$,
  \begin{align}
    \{ a_{n} \langle n! \rangle  \} \equiv 0 \ ({\rm mod} \ {n+1}),
  \label{prime:factorial}
  \end{align}
if $n+1$ is a prime number.
\end{theorem}
For example, $ \{ a_6 \langle n! \rangle  \}  =  9366 = 1338 \times 7 $ and $  \{ a_8 \langle n! \rangle  \}  =  1091670 = 121297 \times 9 - 3$. To prove the above theorem, it is enough to see the following formula holds for a prime number $n$ and a natural number $k$ that is less than or equal to $n-1$,  
  \begin{align}
    k! \times S_{n-1, k} \equiv - (-1)^k \ {\rm mod} \ n,  
  \label{prime:factorial:2}
  \end{align}
because the above formula leads to $\sum_{k=1}^{n-1}  k! \times S_{n-1, k} \equiv 0 \ {\rm mod} \ n$. By using the following formula,
  \begin{align}
    k! \times S_{n-1, k}
      = \sum_{i=0}^k (-1)^i {}_k C_i (k-i)^{n-1},
  \end{align}
and the Fermat's little theorem that states $l^{n-1} \equiv 1 \ {\rm mod} \ n$ if $l$ is not divisible by a prime number $n$, one can show Eq.~(\ref{prime:factorial:2}). Figure \ref{Factorial_prime} shows the example of Eq.~(\ref{prime:factorial:2}) for $n=7$.

We remark that Eq.~(\ref{prime:factorial}) holds for some odd composite numbers. For example, $\{ a_{24} \langle n! \rangle  \} $ can be divided by 25. Also, for all the even numbers, Eq.~(\ref{prime:factorial}) does not hold because $\sum_{k=2}^{n-1}  k! \times S_{n-1, k} \ {\rm mod} \ n$ is an even number mod $n$ and does not equal to $-1 \ {\rm mod} \ n $.

\subsection{Bell number}
Bell number is a number counting how many number of ways to partition of a set of $n$ elements. For example, $\{ a, b, c\}$ can be partitioned into $\{ a \},\{ b \},\{ c \}$ or $\{a\},\{b,c\}$ or $\{b\}, \{c,a\}$ or $\{c\},\{a,b\}$ or $\{a,b,c\}$. Then, the third term of the Bell numbers $B_3$ is 5. The $n$th term of Bell numbers is given by
  \begin{align}
    B_n = \sum_{k=0}^{n-1} {}_nC_k B_k.
  \end{align}
Then, one has
  \begin{align}
    \{ B_n \} = 1, 1, 2, 5, 15, 52, 203, 877, 4140, 21147,  \cdots. 
  \end{align}
We remark that the factorials can be considered as dual of the Bell numbers in terms of the Stirling number. First, the Bell numbers can be represented by using the Stirling number of the second kind as follows.
  \begin{align}
    B_n = \sum_{k=0}^n S_{n,k}.
  \end{align}
Also, the unsigned Stirling number of the first kind $|s(n,k)|$ is a number that counts the number of permutations of $n$ elements with $k$ disjoint cycles. For example, $|s(3,2)| = 3$ because there are three permutations of three elements that fixed the one element: (1,2) or $(2,3)$ or $(3,1)$. Then, the factorial can be represented as follows,
  \begin{align}
    n! = \sum_{k=0}^n |s_{n,k}|.
  \end{align}

Now, we construct the sequence dual of Bell numbers. First, the exponential generating function of the $n$th Bell number is given by
  \begin{align}
    \sum_{n=0}^{\infty} \frac{B_n}{n!} x^n = e^{e^x-1}. 
  \end{align}
Then, by differentiating the each side by $x$, one has
  \begin{align}
    \sum_{n=0}^{\infty} \frac{B_{n+1}}{n!} x^n = e^x e^{e^x-1}.
  \end{align}
Note that $e^x$ can be obtained by substituting $\alpha=1$ into $e^{\alpha x}$. Then, we define the sequence dual of Bell numbers by using $a_1 [e^{\alpha x} ]_L$ that can be obtained by substituting $n=1$ into $\{a_n [e^{\alpha x} ]\}_L$.
  \begin{align}
    \sum_{n=0}^{\infty} \frac{a_n \langle B_n \rangle}{n!} \alpha^n &= \frac{1}{1-\alpha}e^{\left( \frac{1}{1-\alpha} - 1 \right)} =  a_1 [e^{\alpha x} ]_L \, e^{\left( a_1 [e^{\alpha x} ]_L - 1 \right)}.    
  \end{align}
We remark that $\{ a_n \langle B_n \rangle \}$ relates to $B_n$ through the Stirling transform of the first kind corresponding to Eq.~(\ref{factorial_Stirling}) \cite{Stirling:transform}.
\begin{proposition} [Sequence dual of Bell numbers and Stirling transform]
  \begin{align}
    \{ a_n \langle B_n \rangle \} &= \left\{ \sum_{k=0}^n |s(n,k)| \mathcal{S}_1 B_k \right\},
  \\
    &=	1, 2, 7, 34, 209, 1546, 13327, 130922, 1441729, 17572114, \cdots.
  \end{align}
\end{proposition}

We will now focus on the following modular properties of Bell numbers concerning prime numbers. 
  \begin{align}
    B_{n+p} - B_n \equiv B_{n+1} \ \ ({\rm mod} \ p). 
  \end{align}
Then, we have the following modular property for the sequence dual of Bell numbers that should be regarded as sequence dual of the above theorem.
\begin{theorem} [Sequence dual of modular property of Bell numbers]
For natural numbers $n$ and $m$,
  \begin{align}
    \{ a_{n+m}  \langle B_n \rangle \} - \{ a_{n} \langle B_n \rangle \} 
    \equiv 0 \ \  ({\rm mod} \ m).
  \label{seq:Bell:prime}
  \end{align}
\end{theorem}
For example, $a_{7} \langle B_n \rangle - a_{6} \langle B_n \rangle = 130922 - 13327 = 1\times117595$, $a_{7} \langle B_n \rangle - a_{5} \langle B_n \rangle = 130922 - 1546  = 2\times 64688$, $a_{7} \langle B_n \rangle - a_{4} \langle B_n \rangle = 130922 - 209  = 3\times 43571$, $a_{7} \langle B_n \rangle - a_{3} \langle B_n \rangle = 130922 - 34  = 4\times 32722$, $a_{7} \langle B_n \rangle - a_{2} \langle B_n \rangle = 130922 - 7  = 5\times 26183$, $a_{7} \langle B_n \rangle - a_{1} \langle B_n \rangle = 130922 - 2  = 6\times 21820$, $a_{7} \langle B_n \rangle - a_{0} \langle B_n \rangle = 130922 - 1  =  7\times 18703$.

The above theorem can be proved by using the recursion formula,
  \begin{align}
    a_n  \langle B_n \rangle = 2n a_{n-1}  \langle B_n \rangle  - (n-1)^2  a_{n-2} \langle B_n \rangle ,
  \end{align}
and induction method as follows.
  \begin{align}
    a_{n+1+m} \langle B_n \rangle &\equiv 2(n+1) a_n \langle B_n \rangle - n^2 a_{n-1} \langle B_n \rangle \ {\rm mod} \ m \equiv a_{n+1} \langle B_n \rangle \ {\rm mod} \ m.
  \end{align}
Note that for the first two terms of $\{ a_n  \langle B_n \rangle \}$, the following relations hold.
  \begin{align}
    &a_m \langle B_n \rangle \equiv a_0 \langle B_n \rangle \ {\rm mod} \ m \equiv 1 \ {\rm mod} \ m,
  \\
    &a_{m+1} \langle B_n \rangle - a_1 \langle B_n \rangle \equiv a_1 \langle B_n \rangle \ {\rm mod} \ m \equiv 2 \ {\rm mod} \ m.
  \end{align}
The first relation can be proved by the following formula for $\{ a_n \langle B_n \rangle \}$.
  \begin{align}
    a_n \langle B_n \rangle = \sum_{k=0}^n k ! ( {}_n C_k  )^2.
  \end{align}
The second formula can be shown as follows.
  \begin{align}
    a_{1+m} \langle B_n \rangle &= 2(1+m) a_m \langle B_n \rangle - (1+m-1)^2 a_{m-1} \langle B_n \rangle
      \equiv 2 (1+m) \ {\rm mod} \ m
  \nonumber\\
    &\equiv 2 \ {\rm mod} \ m. 
  \end{align}

\newpage
\section{Discussion}\label{Sec:discussion}
In this paper, we created a sequence version of calculus. First, we defined equivalence, some fundamental operations, left/right differential and integral for sequences. Then, we proposed sequence version of power function, exponential function, hyperbolic function, and trigonometric function and sequence versions of Maclaurin series for them. We found the sequence versions of Maclaurin series for exponential function involve divergent series including Grandi's series. We remark that generally there are some candidates of the sequence version of a given function. As for sequence versions of trigonometric functions, first, two kinds of sequence versions of trigonometric functions were proposed, depending on whether the left or right differential was used. The right/left trigonometric functions are graphically represented by the divergent/convergent spiral while the usual trigonometric functions are graphically represented by a circle. Then, we found the sequence version of binomial theorem and Euler's identity. The sequence version of Euler's identity is composed of two equations determined by fourth term of left/right sequence version of $e^{ix}$. Then, by combining the two sequence versions of trigonometric functions, we proposed another sequence version of trigonometric functions whose periodicity is 8. Furthermore, we proposed sequence version of Maclaurin series for Fibonacci sequence, $(P,\ Q)$-Fibonacci sequence, and $k$-bonaaci sequence, and found similarities between them and the sequence versions of exponential function. In addition, we introduce sequence dual of factorial and Bell numbers. Then, the sequence dual of the modular property of factorial concerning the prime numbers (Wilson's theorem) and of the Bell numbers concerning the prime numbers were obtained. Finally, we remark a few points for future research. First, probably, our theory relates systematically to the figurate number including (centered) polygonal numbers, (centered) pyramidal numbers, or (centered) Platonic numbers. Note that $\{ a_n [x^2/2] \}_R$ is the triangle number. Also, there are a few functions that we could not deal with in this paper such as logarithmic function, inverse trigonometric/hyperbolic function, elliptic integral, zeta function. Sequence version of these functions may lead to new formulae for some sequences. Also, considering other frameworks of analogy between sequence and analytics may lead to many new theorems of sequences and analytics.

\appendix
\newpage
\section{List of sequence version of functions and sequence duals}\label{appendix}
Constant function: 
\\
\qquad $\{ a_n [a] \} = a, a, a, a, a, a, a, a, a, a, a, a, \cdots$.
\\
Power function: 
\\
\qquad $\{ a_n [x^k] \} = 0^k, 1^k, 2^k, 3^k, 4^k, 5^k, 6^k, 7^k, 8^k, 9^k, 10^k, 11^k, \cdots$.
\\
Exponential function (right):
\\
\qquad $\{ a_n [e^{\alpha x}] \}_R = (1+\alpha)^0, (1+\alpha)^1, (1+\alpha)^2, (1+\alpha)^3, (1+\alpha)^4, (1+\alpha)^5, (1+\alpha)^6,
\nonumber\\
\qquad\qquad  (1+\alpha)^7, (1+\alpha)^8, (1+\alpha)^9,(1+\alpha)^{10}, (1+\alpha)^{11}, \cdots$.
\\
Exponential function (left): 
\\
\qquad $\{ a_n [e^{\alpha x}] \}_L = (1-\alpha)^0, (1-\alpha)^{-1}, (1-\alpha)^{-2}, (1-\alpha)^{-3}, (1-\alpha)^{-4}, (1-\alpha)^{-5}, (1-\alpha)^{-6},
\nonumber\\
\qquad\qquad (1-\alpha)^{-7}, (1-\alpha)^{-8}, (1-\alpha)^{-9},(1-\alpha)^{-10}, (1-\alpha)^{-11}, \cdots$.
\\
Exponential function (natural): 
\\
\qquad $\{ a_n [e^{-x}] \}^{\rm nat} = 1, 0, 0, 0, 0, 0, 0, 0, 0, 0, 0, 0, \cdots$.
\\
Hyperbolic cosine function:
\\
\qquad $ \{ a_n [\cosh x] \} = 1, \frac{5}{4}, \frac{17}{8}, \frac{65}{16}, \frac{257}{32}, \frac{1025}{64}, \frac{4097}{128}, \frac{16385}{256}, \frac{65537}{512}, \frac{262145}{1024}, \frac{1048577}{2048}, \frac{4194305}{4096}, \cdots$.
\\
Hyperbolic sine function:
\\
\qquad $ \{ a_n [\sinh x] \} = 1, \frac{3}{4}, \frac{15}{8}, \frac{63}{16}, \frac{255}{32}, \frac{1023}{64}, \frac{4095}{128}, \frac{16383}{256}, \frac{65535}{512}, \frac{262143}{1024}, \frac{1048575}{2048}, \frac{4194303}{4096},\cdots$.
\\
Hyperbolic cosine function (natural):
\\
\qquad $ \{ a_n [\cosh x] \}^{\rm nat} = 1, 1, 2, 4, 8, 16, 32, 64, 128, 256, 512, 1024, \cdots$.
\\
Hyperbolic sine function (natural):
\\
\qquad $ \{ a_n [\sinh x] \}^{\rm nat} = 0, 1, 2, 4, 8, 16, 32, 64, 128, 256, 512, 1024, \cdots$.
\\
Cosine function (right):
\\
\qquad $ \{ a_n [\cos x] \}_R  = 1, 1, 0, -2, -4, -4, 0, 8, 16, 16, 0, -32, \cdots$.
\\
Sine function (right):
\\
\qquad $ \{ a_n [\sin x] \}_R  = 0, 1, 2, 2, 0, -4, -8, -8, 0, 16, 32, 32, \cdots$.
\\
Tangent function (right):
\\
\qquad $ \{ a_n [\tan x] \}_R = 0, 1, \infty, -1, 0, 1,-\infty, -1, 0, 1, \infty, -1, \cdots$.
\\
Cosine function (left):
\\
\qquad $ \{ a_n [\cos x] \}_L  = 1, \frac{1}{2}, 0, -\frac{1}{4}, -\frac{1}{4}, -\frac{1}{8}, 0, \frac{1}{16}, \frac{1}{16}, \frac{1}{32}, 0, -\frac{1}{64} \cdots$.
\\
Sine function (left):
\\
\qquad $ \{ a_n [\sin x] \}_L  = 0, \frac{1}{2}, \frac{1}{2}, \frac{1}{4}, 0, -\frac{1}{8}, -\frac{1}{8}, -\frac{1}{16}, 0, \frac{1}{32}, \frac{1}{32}, \frac{1}{64}, \cdots$.
\\
Tangent function (left):
\\
\qquad $ \{ a_n [\tan x] \}_L = 0, 1, \infty, -1, 0, 1,-\infty, -1, 0, 1, \infty, -1, \cdots$.
\\
Cosine function (periodic)
\\
\qquad   $ \{ a_n [\cos x] \}^{\rm per} = 1, \frac{1}{\sqrt{2}}, 0, -\frac{1}{\sqrt{2}}, -1, -\frac{1}{\sqrt{2}}, 0, \frac{1}{\sqrt{2}}, 1, \frac{1}{\sqrt{2}}, 0, -\frac{1}{\sqrt{2}}, \cdots$.
\\
Sine function (periodic)
\\
\qquad  $ \{ a_n [\sin x] \}^{\rm per} = 0, \frac{1}{\sqrt{2}}, 1, \frac{1}{\sqrt{2}}, 0, -\frac{1}{\sqrt{2}}, -1, -\frac{1}{\sqrt{2}}, 0, \frac{1}{\sqrt{2}}, 1, \frac{1}{\sqrt{2}}, \cdots$.
\\
Factorial:
\\
\qquad $ \{ a_n \langle n! \rangle \} = 2, 2, 6, 26, 150, 1082, 9366, 94586, 1091670, 14174522, 204495126, 
\nonumber\\
\qquad\qquad 3245265146, \cdots$.
\\
Bell number:
\\
\qquad $ \{ a_n \langle B_n \rangle \} = 1, 2, 7, 34, 209, 1546, 13327, 130922, 1441729, 17572114, 234662231, 
\nonumber\\
\qquad\qquad 3405357682, \cdots$.



\newpage
\begin{thebibliography}{99}
\bibitem{Fibonacci:history} T. C. Scott and P. Marketos, ``On the origin of the Fibonacci Sequence", MacTutor History of Mathematics, 1 (2014).


\bibitem{Fibonacci:history2} J. A. Adam, Mathematics in nature: Modeling patterns in the natural world, Princeton, NJ: Princeton University Press (2013).


\bibitem{Stirling} L. Comtet, Advanced Combinatorics: The Art of Finite and Infinite Expansions, Dordrecht-Holland{\slash}Boston-U.S.A.: Reidel Publishing Company (1974).


\bibitem{Bell} S. M. Tanny, ``On some numbers related to the Bell numbers", Canadian Mathematical Bulletin ${\bf 17}$, 733 (1975).


\bibitem{Bell2} M. Aigner, ``A characterization of the Bell numbers", Discrete mathematics ${\bf 205}$, 207 (1999).


\bibitem{Bell3} M. Z. Spivey, ``A generalized recurrence for Bell numbers", J. Integer Seq ${\bf 11}$ (2008).


\bibitem{Wilson} G. A. Miller, ``A new proof of the generalized Wilson's Theorem", The Annals of Mathematics ${\bf 4}$, 188 (1903).


\bibitem{Wilson2} S. M. Ruiz, ``An algebraic identity leading to Wilson's theorem", The Mathematical Gazette ${\bf 80}$, 579 (1996).


\bibitem{Wilson3} K. Friedl, and L. R\'{o}nyai, ``Order shattering and Wilson's theorem", Discrete Mathematics ${\bf 270}$, 127 (2003).


\bibitem{Bell:prime} H. W. Becker and J. Riordan, ``The arithmetic of Bell and Stirling numbers", American journal of Mathematics ${\bf 70}$, 385 (1948). 


\bibitem{calculus} M. Abramowitz, I. A. Stegun, Handbook of Mathematical Functions with Formulas, Graphs, and Mathematical Tables, New York: Dover Publications, Ninth printing (1970).


\bibitem{Pascal:binomial} J. L. Coolidge, ``The Story of the Binomial Theorem", The American Mathematical Monthly ${\bf 56}$, 3 (1949).


\bibitem{exponential:infection} S. S. Musa, S. Zhao,  M. H. Wang, A. G. Habib, U. T. Mustapha, and D. He, ``Estimation of exponential growth rate and basic reproduction number of the coronavirus disease 2019 (COVID-19) in Africa", Infect Dis Poverty ${\bf 9}$, 1 (2020).


\bibitem{hyperbolic:relativity} E. M. Lifshitz and L. D. Landau, The Classical Theory of Fields: Volume 2 (Course of Theoretical Physics Series), translated by Morton Hamermesh, fourth revised english edition, Oxford: Butterworth-Heinemann (1980).


\bibitem{trigonometric:rotation} A. Baker, Matrix groups : an introduction to Lie group theory, New York: Springer (2002).


\bibitem{Euler:identity} R. Coolman, Euler's Identity: `The Most Beautiful Equation', livescience (2015).


\bibitem{gamma} E. Artin, The gamma function, New York: Courier Dover Publications (2015).


\bibitem{Diff:Int:seq} D. R. Chalice,``How to Differentiate and Integrate Sequences", The American Mathematical Monthly ${\bf 108}$, 911 (2001).


\bibitem{Yuki:seq} H. Yuki, Search for discrete version of function, Retrieved March 13, 2022, from https:{\slash\slash}www.hyuki.com{\slash}story{\slash}diffsum2.pdf.


\bibitem{Euler:number} L. Carlitz, ``Eulerian numbers and polynomials", Mathematics Magazine ${\bf 32}$, 247 (1959).


\bibitem{Euler:number2} T. Kim, ``Euler numbers and polynomials associated with zeta functions", Abstract and Applied Analysis, Vol. 2008, Hindawi (2008).


\bibitem{Euler:number3} T. K. Petersen, ``Eulerian numbers", Birkh\"{a}user, New York, NY, 3 (2015).


\bibitem{divergent} G. H. Hardy, Divergent series, Vol. 334, American Mathematical Soc (2000).


\bibitem{Grandi} M. Kline, ``Euler and infinite series", Mathematics Magazine, ${\bf 56}$, 207 (1983).


\bibitem{negaFibonacci} J. Triana, ``Negafibonacci Numbers via Matrices", Bulletin of TICMI ${\bf 23}$, 19 (2019). 


\bibitem{Cassini} M. Werman, and D. Zeilberger. ``A bijective proof of Cassini's fibonacci identity", Discrete mathematics ${\bf 58}$, 109 (1986).


\bibitem{gene:Fibo} A. Suvarnamani and M. Tatong, ``Some properties of (p, q)-Fibonacci numbers", Progress in Applied Science and Technology ${\bf 5}$, 17 (2015).


\bibitem{gene:Fibo2} A. Suvarnamani, and M. Tatong, ``Some properties of the product of (P,Q) - Fibonacci and (P,Q) - Lucas number", GEOMATE Journal ${\bf 13}$, 16 (2017).


\bibitem{Pell} W. L. McDaniel, ``Triangular Numbers in the Pell Sequence", Fib. Quart. ${\bf 34}$, 105 (1996).


\bibitem{Pell2} S. Falcon and \'{A}. Plaza, ``The k-Fibonacci sequence and the Pascal 2-triangle." Chaos, Solitons \& Fractals ${\bf 33}$, 38 (2007).


\bibitem{Pell3} C. Bolat and H. K\"{o}se, ``On the properties of k-Fibonacci numbers." Int. J. Contemp. Math. Sciences ${\bf 5}$, 1097 (2010).


\bibitem{Pell4} P. Catarino, ``On some identities for k-Fibonacci sequence." Int. J. Contemp. Math. Sci ${\bf 9}$, 37 (2014).


\bibitem{Jacobsthal} G. E. Bergum, L. Bennett, A. F. Horadam, and S. D.Moore, ``Jacobsthal Polynomials and a Conjecture Concerning Fibonacci-Like Matrices", Fib. Quart ${\bf 23}$, 240 (1985).


\bibitem{Jacobsthal2} A. F. Horadam, ``Jacobsthal representation numbers", significance ${\bf 2}$, 2 (1996).


\bibitem{Jacobsthal3} F. Koken and D. Bozkurt, ``On the Jacobsthal-Lucas numbers by matrix methods", Int. J. Contemp. Math. Sciences ${\bf 3}$, 1629 (2008).


\bibitem{k-bonacci} J. L. Ram\'{i}rez and V. F. Sirvent,``A Generalization of the k-Bonacci Sequence from Riordan Arrays", The electronic journal of combinatorics (2015).


\bibitem{Stirling:transform} M. Bernstein and N. J. A. Sloane, ``Some canonical sequences of integers", Linear algebra and its applications ${\bf 226}$-${\bf 228}$, 57 (1995). 


\end{thebibliography}
\end{document}